\documentstyle[a4wide,amssymb,oldlfont,11pt,epsf]{article}
\textheight23truecm \topmargin-.1truecm \input xy \xyoption{all}

\font\tencmmib=cmmib10 \font\sevencmmib=cmmib7 \font\fivecmmib=cmmib5
\newfam\cmmibfam \textfont\cmmibfam=\tencmmib
\scriptfont\cmmibfam=\sevencmmib \scriptscriptfont\cmmibfam=\fivecmmib
 
\font\script=eusm10 at 12 pt

\def\bs{\backslash}  \def\II{{\bf II}}

 \def\cL{{\cal L}} \def\cO{{\cal O}}
\def\~#1{\tilde{#1}}

\def\gD{\Delta}

\def\gG{{\Gamma}}
\def\gg{\gamma}
\def\gs{\sigma}
\def\gS{\Sigma}

\def\ga{\alpha}
\def\gb{\beta}
\def\gk{\kappa}
\def\gl{\lambda}

\def\go{\omega}

\def\gr{\rho}
\def\grr{\varrho}

\def\bgl{\hbox{\boldmath$\lambda$}}\def\bgm{\hbox{\boldmath$\mu$}}
\def\boldmu{{\bgm}}
\def\bgls{\hbox{\scriptsize\boldmath$\lambda$}}

\def\go{\omega} \def\cO{{\cal O}} \def\gG{\Gamma} 
\def\fZ{{\bf Z}} \def\fR{{\bf R}} \def\fC{{\bf C}} \def\fP{{\bf P}}
\def\cE{{\cal E}} \def\cF{{\cal F}} \def\fF{{\bf F}}
 \def\cX{{\cal X}} \def\cM{{\cal M}}

\def\scM{\hbox{\script M}}

\newtheorem{theorem}{Theorem}[section]

\newtheorem{lemma}[theorem]{Lemma}

\newtheorem{corollary}[theorem]{Corollary}

\def\deg{\hbox{deg}}  
\def\Aut{\hbox{Aut}}  
 \def\det{\hbox{det}} \def\Ker{\hbox{Ker}}

\def\lra{\longrightarrow} 
 
\def\llra{\longleftrightarrow}  \def\Ra{\Rightarrow}
\def\La{\Leftarrow} \def\bs{\ifmmode {\setminus} \else$\bs$\fi}
 \def\tilde{\widetilde} 
\def\inn{\subset}

\def\ende{\hfill $\Box$ \vskip0.25cm } \def\cD{\ifmmode {\cal D}
  \else$\cD$\fi}

\def\xg{\ifmmode {X_{\Gamma}} \else$\xg$\fi} \def\xgeq{\ifmmode {\xg =
    \Gamma\bs\cD} \else$\xgeq$\fi}

\def\xgs{\ifmmode {X_{\Gamma}^*} \else$\xgs$\fi}
           
\def\xgc{\ifmmode {\overline{X}_{\Gamma}} \else$\xgc$\fi}

\begin{document}%

\hfill {BONN--TH--99--02}
\vskip 1cm

\begin{center}
{\Large K3-fibered Calabi-Yau threefolds I,}

{\Large  the twist map} 

\bigskip   
Bruce Hunt 

 {\it and} 

 Rolf Schimmrigk  
\end{center}

\bigskip
\begin{center}
{\Large \bf Introduction}
\end{center}

\bigskip
A natural situation occuring in the general classification theory of
algebraic varieties is that of algebraic fiber spaces $X\lra Y$ such that the
general fiber $F$ has Kodaira dimension 0. This will be the case in
particular if $F$ is Calabi-Yau (a regular variety with trivial canonical
bundle), and if this is the case it is an interesting problem to find
conditions characterizing the case in which the total space $X$ is itself
Calabi-Yau. In this note we make a few remarks on this problem for
three-dimensional $X$. 

Much of what we do is in the context of weighted hypersurfaces. Let $X$ be
a hypersurface of weight $d$ in $\fP_{(w_0,\ldots,w_n)}$. Then a {\it
  sufficient} condition for $X$ to be Calabi-Yau is $d=\sum w_i$. We show
by example (due to I. Dolgachev) that this condition is {\it not
  necessary}. Given a variety $X$ with a fibration $\pi:X\lra Y$, a {\it
  necessary} condition for $X$ to be Calabi-Yau is that $\pi_*K_{X|Y} =
K_Y^{-1}$, and we show by example (of a K3-fibered threefold) that this
condition is {\it not sufficient}. 
Both examples 
 are based on a
particular K3 surface, which is the surface with the largest automorphism
group which preserves the Picard lattice, namely $\fZ/66\fZ$. 

Our interest in this paper is in particular with fibrations whose fibers
have constant modulus. Such varieties are covered by products, and if both
factors of that product are weighted hypersurfaces, this covering can be
neatly described in terms of weighted projective spaces; this is our {\it
  twist map}. This map arose in \cite{HS} in the context of dualities of
heterotic and type II superstrings, but its mathematical formulation gives
a very convenient method for resolving singularities of such
quotients, see Corollary \ref{cquotient}. 
Even in the case of conventional projective spaces, the twist
map is interesting. For example, one corollary is the following perhaps
already known but somewhat startling result.

\medskip
\noindent{\bf Theorem 1:} {\it Let $f(x_1,\ldots, x_n) =x_1^d+\cdots +x_n^d,\
n>3$ be a Fermat polynomial of degree $d$ in $n$ variables, and let
$\{n_1,\ldots, n_{\gl}\}$ be any partition of $n$ with $n_i\geq
2$, $x_{i,k},\ k=1,\ldots,n_{i},\ i=1,\ldots, \lambda$ the
corresponding coordinates. Then the Fermat hypersurface 
$X=\{f=0\}$ is birational to a quotient of
the product $X_1\times \cdots \times X_{\gl}$ of $\gl$ Fermat hypersurfaces
$X_i=\{x_{i,0}^d+x_{i,1}^d+\cdots +x_{i,n_i}^d=0\}$ 
of the same degree $d$ by $\fZ/d\fZ$, acting only on the coordinates
$x_{i,0}$ by multiplication by a primitive $d$th root of unity}.
 
\medskip 
As a particular case of this, consider the Fermat quartic surface $S$ in
$\fP^3$, let $C=\{x_0^4+x_1^4+x_2^4=0\}\inn \fP^2$
(resp. $C'=\{y_0^4-y_1^4-y_2^4=0\}$) be the Fermat quartic
curve in the projective plane. Then under the map 
\begin{eqnarray*} \fP^2\times \fP^2 & \lra & \fP^3 \\
  ((x_0,x_1,x_2),(y_0,y_1,y_2)) & \mapsto & (y_0x_1,y_0x_2,x_0y_1,x_0y_2)
\end{eqnarray*}
the product $C\times C'$ maps onto the surface $S$, displaying the latter as
a $\fZ/4\fZ$-quotient of the former. Recall also that this particular K3
surface is birational to the Kummer surface of the product of elliptic
curves $E_i\times E_{2i}$, where $E_i$ is the elliptic curve with modulus
$\tau = i$ (the unique elliptic curve with an automorphism of order 4) and
$E_{2i}$ is the elliptic curve with modulus $\tau=2i$ (see \cite{In},
p. 546; it is interesting to note that Inose's construction is a
particularly simple case of our twist map). 
On the other hand, there is a natural morphism $\fP^2\lra \fP_{(1,1,2)}$
taking the curve $C$ onto the elliptic curve $E_i$, which is the weighted
hypersurface $y_0^4+y_1^4+y_2^2=0  \inn \fP_{(1,1,2)}$. It is then clear
from construction that the $4-1$ cover of $S$ by $C\times C'$ factors:
\[ \xymatrix{ C\times C' \ar[rr] \ar[dr] & & S \\ & E_i\times E_{2i}\ar[ur]. } \]

Other applications of our twist map are the construction of interesting
examples. In particular, we prove the following {\it existence}
 theorems, the first of which has no analog for elliptic
fibrations or fibrations of abelian surfaces.

\medskip\noindent{\bf Theorem 2:} {\it There exist Calabi-Yau threefolds
  (an example of which is given by 
  a hypersurface of degree 12 in the weighted projective space
  $\fP_{(1,1,2,2,6)}$) which have two different, constant modulus
  K3-fibrations.} 

\medskip\noindent{\bf Theorem 3:} {\it There are examples of Calabi-Yau
  threefolds with both extreme Euler-Poincar\'e characteristics $960$ and
  $-960$, which are images of the twist map, hence have a constant modulus
  fibration. }

\smallskip\noindent
The first of these results Lemma \ref{theorem2}, the second is 
discussed in section
\ref{selliptic} in the text.

\medskip\noindent{\bf Aknowledgements:} The first named author thanks Bernd
Kreussler heartily for his suggesting the general form of the twist map
given here. It is also a pleasure to thank I. Dolgachev and M. Gross for
correspondence in these matters. 

\section{General properties of fibered Calabi-Yau threefolds}

\subsection{Fiber spaces } 
Let $\pi:X\lra Y$ be an algebraic fiber space, i.e., a
proper morphism of algebraic varieties with connected fibers; in general
we will be assuming the base $Y$ of the fibration is smooth and the fibers
are generically smooth. The {\it discriminant} $\gD\inn Y$ is the locus of
$y\in Y$ such that the fiber $X_y$ over $y$ is not smooth, in other words
the image of the set of points for which $d\pi$ fails to have maximal
rank. Assume that $X$ is also smooth, and let $K_X,\ K_Y$ denote the
canonical bundles. The {\it relative canonical bundle} is $K_{X|Y}$,
defined as 
\begin{equation}\label{1} K_{X|Y}=\pi^*K_Y^{-1}\otimes K_X.
\end{equation}
It follows from (\ref{1}) that $\pi_*K_{X|Y} =K_Y^{-1}\otimes \pi_*K_X$, so
in particular
\begin{equation}\label{2} K_X=\cO_X \Ra \pi_*(K_X)=\cO_Y \Ra
 \pi_*K_{X|Y} = K_Y^{-1}.
\end{equation}
This formula gives a {\it necessary} condition on a fiber space $\pi:X\lra
Y$ for $X$ to be Calabi-Yau. Note also that equation (\ref{2}) is an
equation for {\it divisors} on the base $Y$ of the fibration. 

\subsection{Fiber spaces with section } 
Let $\pi:X\lra Y$ be an algebraic fiber space, and assume
now that $\pi$ has a {\it section} $\gs:Y\lra X$. Let $\gS=\gs(Y)$; we can
identify $Y$ with the subvariety $\gS$ of $X$. As such $\gS$ has a normal
bundle $N_X\gS$, and we may apply adjuction, yielding (now as an equation
of divisors on $\gS$, written additively)
\begin{equation}\label{3} {K_X}_{|\gS} = K_{\gS} - c_1(N_X\gS),
\end{equation} 
where $c_1$ denotes the first Chern class, an equivalence class of
divisors. Thus the necessary condition (\ref{2}) takes the form:
\begin{equation}\label{4} c_1(N_X\gS)=K_{\gS}.
\end{equation}

\subsection{Elliptic fibrations } 
Assume now that $X\lra Y$ is elliptic with a section
$\gs:Y\lra X$ and let as above $\gS$ denote the image $\gS=\gs(Y)$. It
follows from fundamental results of Nakayama \cite{Na}, 2.1, that there is
a Weierstra\ss\ model $W(\cL,g_2,g_3)$ over $Y$ and a proper birational map
$\mu:X\lra W(\cL,g_2,g_3)$ over $Y$ such that $\mu(\gS)$ is the {\it zero
  section} of the Weierstra\ss\ model. A Weierstra\ss\ model over $Y$
consists of (a) a line bundle $\cL$ over $Y$, and (b) sections $g_2\in
{\rm H}^0(Y,\cL^{\otimes -4}),\ g_3\in {\rm H}^0(Y,\cL^{\otimes -6})$. The discriminant
$\gD:= g_2^3-27 g_3^2 \in {\rm H}^0(Y,\cL^{\otimes -12})$ is the discriminant of
the Weierstra\ss\ model $W(\cL, g_2,g_3)$, which is defined as follows. Fix
meromorphic sections $x, y, z$ of $\cO_{\fP}(1)\otimes \cL^{-2},\
\cO_{\fP}(1)\otimes \cL^{-3}$ and $\cO_{\fP}(1)$, respectively, where
$\fP:=Proj(\cO_Y\oplus \cL^2\oplus \cL^3)$. Then $W(\cL,g_2,g_3)$ is the
divisor on $\fP$ defined by 
\begin{equation}\label{5} y^2 z = 4x^3 +g_2 x z^2 + g_3 z^3.
\end{equation}
The {\it zero section} of this Weierstra\ss\ model is determined by
dehomoginizing, i.e., choosing the inflection point at infinity of the
Weierstra\ss\ cubic as the zero point of the curve. Note that, if $\gS$
denotes the image of a section as above, then $\cL \cong N_X\gS$, where we
view $\cL$ as a line bundle on $\gS$. Then it follows that $\gD = -12
c_1(N_X\gS)$. 

A Weierstra\ss\ model
is {\it minimal}, if there is no prime divisor $D$ on $Y$ such that
$div(g_2)\geq 4D$ and $div(g_3)\geq 6D$. Nakayama shows that if $Y$ is
smooth and the discriminant $\gD$ has normal crossings, then $W$ has only
rational singularities, if and only if, $W(\cL,g_2,g_3)$ is
minimal.  
\begin{lemma} Let $X\lra Y$ be elliptic with section. The condition
  (\ref{4}) is necessary and sufficient for $X$ to be Calabi-Yau.
\end{lemma}
{\bf Proof:} We must show that $K_{X|\gS}=\cO_{\gS}$ implies
$K_X=\cO_X$. This follows from Kodaira's formula for the canonical bundle,
in the higher-dimensional formulation as in \cite{Ka}, which is\footnote{in
  general there will be an additional error term which is present due to
  less accurate control over the birational geometry of these spaces in
  higher dimensions; by results of Grassi \cite{G} this term can be avoided
  by choosing our model correctly}
\begin{equation}\label{6} K_X = \pi^*(K_Y + \sum a_i[\gD_i]),
\end{equation}
where $\gD_i$ are the irreducible components of the discriminant $\gD$ and
the rational numbers $a_i$ are determined by the type of singular fiber
over $\gD_i$. For every $a_i$ we have $12 a_i \in \fZ$, and setting $\gD =
\sum 12 a_i \gD_i$, the divisor $\gD$ is divisible by 12 with $\gD = \sum
12 a_i \gD_i = -12 c_1(N_X\gS) = -12c_1(\cL) \Ra c_1(\cL)=-\sum a_i \gD_i$. At
any rate, (\ref{6}), (\ref{4}) and the assumption $K_{X|\gS} = \cO_{\gS}$
together give $K_X = \pi^*(K_{X|\gS}) = \pi^*(\cO_{\gS}) = \cO_X$.  \ende

Finally, if $X\lra Y$ is an elliptic fibration and $X$ is a smooth threefold
with $K_X=\cO_X$, then $X$ is birational to a minimal Weierstra\ss\ model
over a surface $S$ which is one of the following:
\begin{enumerate}\item a minimal surface with $\gk(Y)=0$, $g_2,\ g_3$
  constant. \item a minimal ruled surface over an elliptic curve. \item
  $S=\fP^2$ or $S$ is one of the Hirzebruch surfaces $\fF_n,\ 0\leq n \leq 12$.
\end{enumerate}
(For $n> 12$, any Weierstra\ss\ model is necessarily non-minimal, hence by
the above result has non-rational singularities and is consequently not
Calabi-Yau).

\subsection{Abelian surface fibrations } 
Next we assume that the general fiber of $X\lra Y$ is an
abelian surface, and let $\gs:Y\lra X$ be a section. Now the image $\gS$
has codimension two, being just a point in each fiber. The normal bundle
$N_X\gS$ is a rank two bundle on $\gS$ (and on $Y$), and $c_1(N_X\gS)$ is a
divisor on $\gS$, so the formula (\ref{4}) still makes sense. Here we also
have a formula for the canonical bundle \cite{A}, 2.16, which is the
following:
\begin{equation}\label{7} K_X = \pi^*(K_Y + \sum a_i[\gD_i]) +
  \sum_i^N(e_jF_j + E_j) + \sum_1^M (d_iD_i + D_i'),
\end{equation}
where the first sum is over all singular fibers, the second sum over all
{\it multiple} fibers and the third sum is over all components of the
singular fibers which contain {\it two or more} components. It follows from
this that we have 
\begin{lemma} Let $X\lra Y$ be a fiber space of abelian varieties over a
  curve, and suppose that 
  \begin{enumerate} \item There are no multiple fibers. \item All singular
    fibers are irreducible. 
  \end{enumerate}
  Then the condition (\ref{4}) is necessary and sufficient for $X$ to be
  Calabi-Yau. 
\end{lemma}
{\bf Proof:} Using the formula (\ref{7}), the proof is just as
above, the important fact being that under the stated assumptions $K_X$ is
the pullback of a class on the base. \ende 

\subsection{K3-fibrations}
In the case that $X$ is a K3-fibration, assuming $X$ is a Calabi-Yau
threefold, the base $Y$ is necessarily $\fP^1$. The necessary condition
(\ref{4}) becomes in this case $c_1(N_X\gS) =-2$. As we shall see below,
however, this condition is {\it not} sufficient, not even under assumptions
as in the previous lemma. We will give an example of a K3-fibration, for
which (\ref{4}) is satisfied, without multiple fibers, and for which all
singular fibers are irreducible, but which is far from being Calabi-Yau. To
discuss this example, we first discuss monodromy and how this relates to
the coefficients $a_i$ which occur in (\ref{6}) and (\ref{7}). 

\section{Torsion monodromy and fibrations with constant modulus}
\subsection{Monodromy}
Let $X\lra Y$ be an algebraic fiber space, and let $Y_0\inn Y$ denote the
Zariski open subset of $Y$ over which $\pi$ is smooth, i.e., the complement
of the discriminant $\gD$.
Let $\pi_0:X_0\lra Y_0$ denote the
restriction of $\pi$ to $X_0=X-\pi^{-1}(Y-Y_0)$, so that $\pi_0$ {\it is} a
locally trivial fibration in the sense of topology. The collection of the
integral homology groups of the fibers, ${\rm H}^k(F_y,\fZ)$, ($y\in Y_0$) form a
sheaf $
\cF_0$ over $Y_0$, denoted $R^k(\pi_0)_*\fZ$, and the
monodromy around any $s\in \gD:=Y-Y_0$ can be viewed as a translation in
the fiber over some fixed base point $*\in Y_0$, $(\cF_0)_*$; sending
$\gamma \in \pi_1(Y_0,*)$ to the matrix $T_{\gamma}\in \Aut((\cF_0)_*) =
\Aut({\rm H}^k(F_*,\fZ))$ which describes this translation gives the monodromy
representation \[ \rho:\pi_1(Y_0,*) \lra \Aut((\cF_0)_*).\]
By the monodromy theorem, each monodromy matrix around an isolated bad
fiber has eigenvalues which are roots of unity. Let $T$ be the monodromy
matrix; if $(T-1)$ is {\it nilpotent}, one says $T$ is {\it
  unipotent}. Suppose $T^k=1$ for some $k$; in this case we will speak of
{\it torsion} monodromy. The general monodromy matrix is {\it
  quasi-unipotent}, meaning that for some $N$ and $k$ we have
$(T^k-1)^N=0$. In this case, $N$ is called the degree of unipotency, and it
is at most the dimension of the fiber. 

Suppose that we have a {\it local} fiber space $\cX \lra D$ over the unit
disc $D$, and assume that all fibers $\cX_z$ are smooth for $z\in D^*=D -
\{ 0 \}$. Assume furthermore that the monodromy is torsion, say $T^k=1$,
and consider the base change $D\lra D,\ z\mapsto z^k$. Pulling back the
fibration $\cX$, we see that the monodromy is now $T^k=1$, i.e.,
trivial. Again, if $X\lra Y$ is an algebraic fiber space over a curve $Y$,
and all monodromies are torsion, then we can find a (branched) cover
$Y'\lra Y$, such that the monodromy of the pull-back of $\cX$ is
trivial. For this it is sufficient to have the branching of order $k_i$ at
each point $z_i$ for which the local monodromy matrix $T_i$ has order
$k_i$. 

While the sheaf $\cF_0$ is a sheaf of $\fZ$-modules, by tensoring with
$\cO_{Y_0}$ we get a sheaf of $\cO_{Y_0}$-modules, call it 
\[\cE_0:=\cF_0\otimes \cO_{Y_0}.\]
There is a very special extension of $\cE_0$, called the {\it canonical
  extension} (Schmidt) $\cE$ of $\cE_0$ to $Y$, i.e., such that
$\cE_{|Y_0}=\cE_0$. It is known how to describe holomorphic sections
generating $\cE$: let $r=rk_{\fZ}{\rm H}^k(F,\fZ)$, and let $v_1,\ldots, v_r$ be
a $\fZ$-basis. For each $s\in \gD$ let $T=T_s$ be the monodromy map
defined by a small loop around $s$, and let $t=t_s$ be a local coordinate on
the base near $s$ (i.e., $t=0$ defines the point $s\in Y$). Then, for
$j=1,\ldots,r$, the expressions 
\[ \sigma_j=\exp\left({1\over 2\pi i}\log T\log t\right) v_j\]
define holomorphic sections of $\cE$ and in fact generate it. If this
expression seems somewhat formidable\footnote{The first named author found
  it so; he is indebted to Donu Arapura for the following explanation}, 
consider the case that we only have
torsion monodromy. 
Consider a Galois cover of $Y$ such that for each $s_i\in \gD$, if $k_i$
denotes the order of $T_{s_i}$ ($T_{s_i}^{k_i}=1$), then the cover
$p:Y'\lra Y$ is given locally by $t_{s'}'\mapsto (t_{s'}')^{k_i}=t_s$ for
any $s'\in p^{-1}(t_{s_i}),\ s_i\in \gD$. By the results above, the lift of
$X\lra Y$ to $Y'$, 
\[\begin{array}{rcccl} & X' & \stackrel{p}{\lra} & X & \\ 
\pi' & \downarrow & & \downarrow & \pi \\
& Y' & \stackrel{p}{\lra}  & Y & 
\end{array}\]
will have trivial monodromy. The eigenvalues of $T$ are of the form
$e^{2\pi i \alpha_j},\ j=1,\ldots, r$, where $\alpha_j=p_j/k_i$ are
rational numbers with $k_i$, $p_j\in \fZ$ (not necessarily in
lowest terms). Let us suppose that $T$ is diagonalized; 
then the matrix $ \log T$ can be calculated as follows:
\[ \log T = \log \left(\begin{array}{ccc} e^{2\pi i \alpha_1} & & 0 \\ &
    \ddots & \\ 0 & & e^{2\pi i \alpha_r}
  \end{array}\right) = \left(\begin{array}{ccc} 2\pi i \alpha_1 & & 0 \\ &
    \ddots & \\ 0 & & 2\pi i \alpha_r
  \end{array}\right),\]
so that the expression  $ \exp\left({1\over 2\pi i}\log T\log
  t\right)$ can be calculated as 
\[  \exp\left({1\over 2\pi i}\log T\log
  t\right) =\left(e^{\log t}\right)^{{1\over 2\pi i}\log T} =
t^{\left({1\over 2\pi i} \log T\right)} = \left(\begin{array}{ccc}
    t^{\alpha_1} & & 0 \\ & \ddots & \\ 0 & & t^{\alpha_r}
  \end{array}\right), \]
and we have $\sigma_j= \left(\begin{array}{ccc}
    t^{\alpha_1} & & 0 \\ & \ddots & \\ 0 & & t^{\alpha_r}
  \end{array}\right) v_j$ for $j=1,\ldots, r$. So in this case we see how
to calculate generators of $\cE$ -- they are just certain (rational) powers
of the local coordinate times the locally constant section. 

Now consider the sheaves $\pi_*K_{X|Y}$ and $\pi'_* K_{X'|Y'}$.
This sheaf is related to the canonical extension $
\cE$ above by the following result of Kawamata (\cite{Ka2}, Theorem 1):
\[ \pi_* K_{X|Y} \cong i_* \cE_0^{k,0} \cap \cE,\]
where $i:Y_0\rightarrow Y$ is the inclusion and $\cE_0^{k,0}\inn \cE_0$ is
the ${\rm H}^{k,0}$-part of the Hodge decomposition of ${\rm H}^k(F,\fC)$. Note this is
an isomorphism of line bundles, and the monodromy group acts here with a
single root of unity at each $\gD_i\inn \gD$; this root of unity
$e^{2\pi i \alpha_i}$ gives the coefficient of $\gD_i$ in the sheaf
$\pi_* K_{X|Y}$ above. 
We now apply this to K3 fibrations $X\lra \fP^1$, the $\gD_i$ are points
on $\fP^1$, $k=2$, and the result above shows that:
\begin{quote} the action of the monodromy on the holomorphic two-form
  determines the coefficient $r_i$ of $\gD_i$ in $\pi_*K_{X|Y}$ by the
  rule:
\begin{equation} \label{formula} \hbox{If } T_i(\omega) = 
  e^{2\pi i \alpha_i} \omega \quad \hbox{( $\omega$ the
  holomorphic two-form)},
\end{equation}
then $r_i = \alpha_i$, i.e., $\pi_*K_{X|Y} = \sum \alpha_i \gD_i$. In
particular, $\sum \alpha_i =2$ is a necessary condition for $X$ to be
Calabi-Yau. 
\end{quote}

\subsection{Constant modulus} In general, given an algebraic fiber space
$X\lra Y$, the complex structure (modulus) of the fibers vary, in a
holomorphic manner. However, it can also occur that the modulus is fixed; in
this case we speak of {\it constant modulus}. Easy examples are given by
elliptic surfaces for which the $J$-invariant is a constant; indeed,
suppose $X\lra Y$ is an elliptic fibration, and let $W\lra Y$ be the
Weierstra\ss\ model. Then $J=g_2^3 / \gD$, and $J$ is constant when
$g_2$ or $g_3$ vanish. For example, suppose $g_2=0$, so that $J\equiv 0$. Then
all singular fibers are of types $II,\ II^*, IV, IV^*$. More precisely,
suppose we choose $Y=\fP^1$, and let $\gD$ consist of two points, with
singular fibers of type $II$ at one and $II^*$ at the other. This is a
rational elliptic surface, and clearly all elliptic fibers have $J =0$ and
so are copies of
the elliptic curve with automorphism group $\fZ/6\fZ$. In this case by the
result of the previous section, there is a cover $Y'\lra \fP^1$, such that
the pull back of the fibration has trivial monodromy. Furthermore, the
modulus is still constant. Under these circumstances, it is clear that the
pull back of the fibration is a fibration without monodromy and with
constant modulus. It follows that this pull back is a {\it product}. This
fact holds more generally.

\subsection{Uniformisation} Let $X\lra Y$ be an algebraic fiber space
satisfying the following conditions:
\begin{enumerate} \item The dimension of $Y$ is one. 
  \item All monodromies around bad fibers are
  torsion. \item The modulus of the fibers of $X$ is constant. 
\end{enumerate}
\begin{lemma} Under the assumptions just made, there is a finite, branched
  cover $Y'\lra Y$ such that the pull back of $X$ to $Y'$ is birational to 
  a product. 
\end{lemma}
{\bf Proof:} For this it is sufficient to construct a cover $Y'\lra Y$,
which is branched at each of the base points $y_i\in Y$ 
of the bad fibers to degree
$k_i$, where $T_i^{k_i}=1$, $T_i$ the local monodromy matrix at the point
$y_i$. Since the dimension of $Y$ is one, this can clearly be done. \ende

\noindent{\bf Remark:} Under some mild assumptions of the moduli space of
the fiber in question (as in, for example, \cite{jraw}), the assumption 3.
{\it implies} the assumption 2.\ above. If this is the case and the general
fiber is smooth, the assumption 3.\ means that the period map has image
a point which is in the {\it interior} of the moduli space. Since the
corresponding moduli point is for a smooth variety, the only degeneration
which can occur is a torsion one, gotten essentially by taking a finite
quotient of a smooth fiber. Thus 3.\ $\Ra$ 2. \ende 

\section{A construction of quotients as weighted projective hypersurfaces}
\subsection{Notations}
We will be working with weighted projective spaces, which are certain
(singular) quotients of usual projective space. Alternatively, they may be
described as quotients of $\fC^{n+1}$ by a $\fC^*$-action. We assume the
weights $(w_0,\ldots, w_n)$ are given, let $\boldmu_{w_i}$ denote the group
of $w_i$th roots of unity, and consider the action of $\boldmu:=
\boldmu_{w_1}\times \cdots \times \boldmu_{w_n}$ on $\fP^n$ as follows. Let
$g=(g_0,\ldots, g_n) \in \boldmu$, and consider for $(z_0:\ldots: z_n)$
homogenous coordinates on $\fP^n$ the action
\[ (g,(z_0:\ldots: z_n)) \mapsto (g_0z_0:\ldots: g_nz_n).\]
Alternatively, consider the action of $\fC^*$ on $\fC^{n+1}$ given by 
\[ (t,(z_0,\ldots, z_n)) \mapsto (t^{w_0}z_0,\ldots, t^{w_n}z_n).\]
In both cases, the resulting quotient is the weighted projective space,
which we will denote by $\fP_{(w_0,\ldots, w_n)}$. General references for
weighted projective spaces are \cite{dolg} and \cite{Y}. A {\it weighted
  hypersurface} is the zero locus of a weighted homogenous polynomial
$p$. Such a hypersurface or the corresponding polynomial is called {\it
  transversal}, if the only singularities are the intersections with the
singular locus of the ambient weighted projective space, and quasismooth,
if the cone over the hypersurface is quasismooth, i.e., smooth outside of
the vertex; for weighted {\it hypersurfaces}, these notions are equivalent
(cf.\ \cite{Y}, Propositions 6 and 8). 

 We will assume the weights are {\it normalized} in the sense that no
$n$ of the $n+1$ weights have a common divisor $>1$. Both for the weighted
projective spaces as well as for the weighted hypersurfaces this assumption
is no restriction (cf.\ \cite{dolg} 1.3.1 and \cite{Y}, pp.\ 185-186). 
We will write such isomorphisms in the sequel without further comment, for
example $\fP_{(2,3,6)}\cong \fP_{(2,1,2)}\cong \fP_{(1,1,1)}=\fP^2$, where
the first equality is because the last two weights are divisible by 3, the
second while the first and last are divisible by 2. 

We will use the notation $\fP_{(w_0,\ldots, w_n)}[d]$ to denote either a
certain weighted hypersurface of degree $d$, or to denote the whole family
of such (the context will make the usage clear). In the particular case
that the weighted polynomial $p$ is of Fermat type, then there is a useful
fact, corresponding to the above normalizations. For example, in
$\fP_{(2,3,6)}$ consider the weighted hypersurface $x_0^6+x_1^4
+x_2^2=0$. Then the isomorphism $\fP_{(2,3,6)}\cong \fP_{(2,1,2)}$ above is
given by the introduction of new variables $(x_0')=x_0^3,$
which is in spite of appearances a one to one coordinate transformation
(becuase of admissible rescalings), and the Fermat polynomial becomes
$(x_0')^2+x_1^4+x_2^2=0$. Again, the isomorphism $\fP_{(2,1,2)}\cong
\fP_{(1,1,1)}$ is given by setting $(x_1')=x_1^2$, and the Fermat
polynomial becomes $(x_0')^2+(x_1')^2+x_2^2=0$, which is a quadric in the
projective plane. We {\it denote} this process by the symbolic expressions
\[ \fP_{(2,3,6)}[12]\cong \fP_{(2,1,2)}[4] \cong \fP_{(1,1,1)}[2].\]

\subsection{The construction}
We now introduce the twist map; this map will give an explicit form to the
forming of quotients of products $V_1\times V_2$ of weighted hypersurfaces
by an abelian group acting on the product. 

\subsubsection{The twist map} Let $V_1,\ V_2$ be weighted hypersurfaces
defined as follows.
\begin{equation}\label{8} \parbox{10cm}{$\displaystyle V_1 = \{x_0^{\ell}
    +p(x_1,\ldots,x_n) =0 \} \inn \fP_{(w_0,w_1,\ldots, w_n)} \\ V_2 = \{
    y_0^{\ell} + q(y_1,\ldots, y_m) =0 \} \inn \fP_{(v_0,v_1,\ldots,v_m)}$},
\end{equation}
where we assume both $p$ and $q$ are quasi-smooth. The degrees of these
hypersurfaces are 
\[ \nu=\deg(V_1) = \ell\cdot  w_0,\quad \mu=\deg(V_2) = \ell\cdot v_0.\]
We then consider the hypersurface 
\begin{equation}\label{9} X:= \{ p(z_1,\ldots, z_n) -q(t_1,\ldots, t_m) =0
  \} \inn \fP_{(v_0 w_1, \ldots, v_0 w_n, w_0 v_1, \ldots, w_0 v_m )}.
\end{equation}
Note that the degree of $X$ is $v_0\cdot \deg(p)=w_0\cdot \deg(q) =
v_0w_0\ell$. 
\begin{lemma}\label{l3.1} The rational map 
  \begin{eqnarray*} \Phi: \fP_{(w_0,w_1,\ldots, w_n)} \times
    \fP_{(v_0,v_1,\ldots,v_m)} & \lra & \fP_{(v_0 w_1, \ldots, v_0 w_n, w_0
    v_1, \ldots, w_0 v_m )} \\ ((x_0,\ldots, x_n),(y_0,\ldots, y_m)) &
    \mapsto & (y_0^{w_1/w_0}\cdot x_1, \ldots, y_0^{w_n/w_0}\cdot x_n,
    x_0^{v_1/v_0}\cdot y_1, \ldots, x_0^{v_m/v_0}\cdot y_m)
  \end{eqnarray*}
restricts to $V_1\times V_2$ to give a rational generically finite
 map onto $X$. 
\end{lemma}
{\bf Proof:} First we show that the map is well-defined outside of the
locus\footnote{As the point $\{(1,0,\ldots, 0), (1,0,\ldots,0)\}$ is not
  contained on $V_1\times V_2$ we will disregard this in what follows}
 \[\{y_0=x_0=0\}\cup \{(1,0,\ldots, 0), (1,0,\ldots,0)\}
 \inn \fP_{\bf w}\times \fP_{\bf v},\] where we have used the abbreviations
 $\fP_{\bf w}$ and $\fP_{\bf v}$ for the projective spaces above; similarly
 we shall use the abreviation $\fP_{\bf w,v}$ for the image projective
 space. That $\Phi$ is well-defined as claimed holds because the variables
 $z_1,\ldots, z_n$ (resp. $t_1,\ldots, t_m$) have weights all divisible by
 $v_0$ (resp. by $w_0$), hence a change of the branch of the $w_0^{th}$
 roots of $y_0$ (resp. of the $v_0^{th}$ roots of $x_0$) just amounts to an
 admissible overall scaling of the coordinates. The locus $\{y_0=x_0=0\}$
 is the locus where all image values are zero, hence this consitutes the
 locus where $\Phi$ is not defined (put differently, where $\Phi$ is not a
 morphism but only a rational map). If we restrict $\Phi$ to $V_1\times
 V_2$, then $x_0^{\ell} = -p(x_1,\ldots, x_n)$ and
 $y_0^{\ell}=-q(y_1,\ldots, y_m)$. Hence
\begin{eqnarray*} p(z_1,\ldots, z_n)-q(t_1,\ldots, t_m) & = &
  p(y_0^{w_1/w_0}\cdot x_1, \ldots, y_0^{w_n/w_0}\cdot x_n) -
  q(x_0^{v_1/v_0}\cdot y_1, \ldots, x_0^{v_m/v_0}\cdot y_m) \\ & = &
  y_0^{\nu/w_0}p(x_1,\ldots, x_n) - x_0^{\mu/v_0}q(y_1,\ldots, y_m) \\ & =
  &  y_0^{\ell}p(x_1,\ldots, x_n) -x_0^{\ell} q(y_1,\ldots, y_m) \\ & = &
  -q(y_1,\ldots, y_m)\cdot p(x_1,\ldots, x_n) +p(x_1,\ldots, x_n)\cdot
  q(y_1,\ldots, y_m) =0.
\end{eqnarray*}
Since $V_1\times V_2$ and $X$ both have dimension $n+m-2$, it is clear that
$\Phi_{|V_1\times V_2}$ is finite-to-one onto its image.  \ende

We call the rational map $\Phi$ the {\it twist map}.

Let $\boldmu_{\ell}$ denote the group of the $\ell$th roots of unity in
$\fC$. 
\begin{corollary} Assume that gcd$(w_0,v_o,\ell)=1$. Then 
$\Phi$ maps $V_1\times V_2$ to the quotient $V_1\times
  V_2 / \boldmu_{\ell}$, where $\boldmu_{\ell}$ acts effectively 
on $ V_1\times V_2 \inn \fP_{(w_0,\bf
  w)} \times \fP_{(v_0,\bf v)}$ by ($\gg \in \boldmu_{\ell}$) 
\[ (\gg,(x_0:\cdots :x_n),(y_0:\cdots:y_m)) \mapsto ((\gg
  x_0:x_1:\cdots:x_n), (\gg y_0:y_1:\cdots:y_m)).\]
Hence the map $V_1\times V_2 \lra X$ is generically $\ell : 1$. 
\end{corollary}
{\bf Proof:} The product $V_1\times V_2$ is given by the two equations $\{
x_0^{\ell} + p(x_1,\ldots, x_n) = 0 = y_0^{\ell} + q(y_1,\ldots, y_m) \}$,
hence it is invariant under the given action. Suppose
$((x_0,\ldots,x_n),(y_0,\ldots, y_m))$ is a solution of the equations; then
for any $\gg \in \boldmu_{\ell}$, $((\gg x_0,\ldots,x_n),(\gg y_0,\ldots,
y_m))$ is also a solution. Provided the weights $w_0,\ v_0$ and $\ell$ have
no common divisor, none of the $\gg$ act as an admissible overall scaling,
hence the result. \ende

\subsubsection{Resolving quotients}
It is well-known that the weighted projective spaces have only quotient
singularities (see \cite{dolg}, 1.2.5 and 1.3.3) 
which can be resolved by the methods of torus
embeddings. Furthermore, from our assumption that $p$ and $q$ are
quasismooth, it follows that also $V$ has only quotient singularities (cf.\
\cite{dolg}, 3.1.6). This is then also true of the polynomial $p-q$
defining $X$, hence $X$ also has only quotient singularities, which can
again be resolved by torus methods. 

We have seen that there is an action of $\boldmu_{\ell}$ on $V$; we let
$\tilde{V}$ be a resolution of $V$ to which the action of $\boldmu_{\ell}$
lifts. Let $\tilde{X}$ be a resolution of $X$. From these assumptions, the
map $\Phi$ lifts to a map $\tilde{\Phi}: \tilde{V} \lra \tilde{X}$, and we have
the following commutative diagram
\[ \xymatrix{ \tilde{V}
 \ar[r]^{\tilde{\Phi}} \ar[d] & \tilde{X}\ar[d] \\ V \ar@{
    -->}[r]^{\Phi} & X.}\] 
By assumption both $\tilde{V}$ and $\tilde{X}$ are
smooth.
\begin{corollary}\label{cquotient}
  Under the above assumptions, $\~X$ is a resolution of the
  quotients $V_1\times V_2/\boldmu_{\ell}$ and $\tilde{V}/\boldmu_{\ell}$. 
\end{corollary}

\subsubsection{The fibration}
If we project the quotient $V_1\times V_2/\boldmu_{\ell}$ onto the
individual factors, we get two {\it rational} fibrations, $\~X  
\xymatrix{\ar@{-->}[r] & V_1/\boldmu_{\ell}} 
$ and $\~X \xymatrix{\ar@{-->}[r]& V_2/\boldmu_{\ell} }$. In each case, the
generic smooth fibers are copies of resolutions of 
 $V_2$ (resp.\ $V_1$). We are especially
interested in the case that $\~X$ is Calabi-Yau.
\begin{lemma}\label{l3.4}
  Suppose $\~X$ is Calabi-Yau and $w_0>1$. Then the rational 
  fibration $\~X \xymatrix{\ar@{-->}[r]&   V_1/\boldmu_{\ell}}$
  induces a genuine fibration onto a resolution $Y$ of
  $V_1/\boldmu_{\ell}$, if and only if $\tilde{V}_2$ is also
  Calabi-Yau. 
\end{lemma}
{\bf Proof:} ``$\Ra$'': if $\~X \lra Y$ is a fibration and $c_1(\~X) = 0$,
then by adjunction $c_1(F)= c_1(\~V_2) = c_1(\~X)_{|F} - c_1(N_{\~X}F) =0$.

``$\La$'': Suppose $c_1(\~X) =c_1(F)=c_1(\~V_2) =0$. Since by adjunction
this implies $c_1(N_{\~X}F)=0$, it suffices to show that $X \lra
V_1/\boldmu_{\ell}$ lifts to a morphism $\~X \lra Y$ for a resolution $Y$
of $V_1/\boldmu_{\ell}$. The map $X\lra V_1/\boldmu_{\ell}$ is given by
fixing $(x_1:\cdots:x_n)$ and mapping all $(y_0:\cdots :y_m)$ with
$\Phi(x,y) \in X$ to $(x_1:\cdots:x_n)$. Thus the projection is
well-defined unless $y_0=0$, a locus which is however blown up upon
resolution of $X$, provided $w_0>1$. \ende

\begin{corollary} \label{cfiber} Let $V_2$ and $X$ fulfill the necessary
  conditions for being Calabi-Yau, $\sum_{j=0}^m v_j = \ell v_0$ and
  $v_0\sum_{i=1}^n w_i + w_0 \sum_{j=1}^m v_j = v_0 w_0 \ell$. Then $X$ has
  a resolution of singularities $\~X$ which is a fiber space over $Y$ with
  constant modulus.
\end{corollary}
We now consider the possible singular fibers which can occur. 
\begin{lemma}\label{l3.9} 
  The singular fibers of $\~X \lra Y$ occur at the (image in
  $V_1/\boldmu_{\ell}$ of the) set of fixed points of the $\boldmu_{\ell}$
  action on $V_1$.
\end{lemma}
{\bf Proof:} This is well-known. \ende
We must consider two situations. First, if $x_0 =0$, then the entire locus
$\~{\gD} = \{ p(x_1,\ldots,x_n)=0\} \inn V_1$ is contained in the
discriminant and maps to $\gD \inn
V_1/\boldmu_{\ell}$, and we let $\gr^*(\gD)$ denote the total transform on
$Y$, where $\gr:Y\lra V_1/\boldmu_{\ell}$ is the resolution induced by that
of $\fP_{(\bf v,\bf w)}$. Secondly, it can happen that for $x_0\neq 0$ we
have further fixed points, a phenomenon which however only occurs in the
case of weighted projective spaces (we will show examples below). It is
easy to see that this can only occur along loci of the type $x_i=0$ for
some $i\in \{1,\ldots,n\}$. 

\noindent{\bf Remark:} {\it Finding the fibers of the fibration}

In the physics literature these fibrations and the fibers are often found
by the method of ``eliminating coordinates''. As this is a source of
confusion, we remark here on this. Sometimes it works, but often one gets
misleading results. Consider the example of the K3 surface given
by the equation
\[ X= \{ z_1^{12}+z_2^6+z_3^4+z_4^2 =0 \}\inn \fP_{(1,2,3,6)}.\]
According to our twist map, this is the quotient of the product 
\[ \{x_0^4+x_1^{12}+x_2^6 =0\}\times \{ y_0^4+y_1^4+y_3^2=0\} \inn
\fP_{(3,1,2)} \times \fP_{(1,1,2)}.\] 
The latter curve is elliptic, and by our results above, the quotient of the
product by $\fZ/4\fZ$ has a fibration with constant fibers equal to that
elliptic curve. 
For $\bgl = (\gl_1:\gl_2)\in \fP_{(1,2)}$ let $D_{\bgls}=\{
\gl_2z_1^2-\gl_1^2 z_2 =0 \} \inn \fP_{(1,2,3,6)}$. Let $X_{\bgls} = X\cap 
D_{\bgls}$. This actually defines the fibration we already derived above,
over the weighted projective line $\fP_{(1,2)}$ (which is of course
isomorphic to the usual $\fP^1$). The equation for the fiber can be derived
by eliminating $z_2$, but upon eliminating $z_1$ (which occurs
quadratically in the equation of $D_{\bgls}$), one gets the equation of a
{\it quotient} of the fiber. These equations are: 
\[X_{\bgls} =\{ (1+\left({\gl_1^2\over
\gl_2}\right)^6)z_2^6+z_3^4+z_4^2=0\}\inn \fP_{(2,3,6)}\cong 
\{(1+\left({\gl_1^2\over \gl_2}\right)^6) 
(z_2')^2+(z_3')^2+z_4^2=0\}\inn  \fP_{(1,1,1)},\]
which describes a rational curve (a quotient of the actual fiber), and 
\[X_{\bgls}=\{((1+\left({\gl_2\over
\gl_1^2}\right)^6)z_1^{12}+z_3^4+z_4^2=0\}\inn
\fP_{(1,3,6)}\cong\{(1+\left({\gl_2\over
\gl_1^2}\right)^6)
(z_1')^4+z_3^4+z_4^2=0\}\inn  \fP_{(1,1,2)},\]
which describes the elliptic curve which is the fiber. The rule to follow
is then: if one of the weights of the eliminated variables is one,
and the variable which is eliminated occurs {\it linearly} in the expression of
the divisor $D_{\bgls}$, then one gets the correct answer. 
To see that one must have one of the weights equal to unity, we consider
one more example. Consider the product of two weighted hypersurfaces of
degrees 6 and 12, respectively, in the product $\fP_{(2,1,1)}\times
\fP_{(4,1,1,6)}$, and let $\boldmu_3$ act as above on this product
$V_1\times V_2$. Clearly $V_2$ is a K3 surface, and we get a fibration over
$\fP_{(1,1)}$, whose fibers are isomorphic to $V_2$. By the twist map this
maps to a hypersurface of degree 24 in $\fP_{(4,4,2,2,12)}$, which is
isomorphic to a hypersurface of degree 12 in 
$\fP_{(2,2,1,1,6)}$. Note that in the equation of the
corresponding divisor $D_{\bgls}$, both coordinates $z_1$ and $z_2$ occur
linearly, but both have weight 2. Upon elimination of either, we get a
hypersurface of degree 12 in $\fP_{(2,1,1,6)}$, which is of general type
and not K3. However, it is easy to see that it is a $2-1$ {\it cover} of
the K3 we started with. Hence the linearity of the coordinate in the
equation of the divisor is not sufficient.

\subsection{K3 surfaces}
We now apply the results above to the construction of K3 surfaces as
quotients of the type $C\times E$, where $C$ is a curve with an action of
$\fZ/\ell\fZ$, and $E$ is an elliptic curve with the same automorphism. To
apply the above, we must present both as weighted hypersurfaces.

\subsubsection{Elliptic curves} We consider the following elliptic curves.
\begin{eqnarray*} E_1 & = & \{ y_0^3 +y_1^3 +y_2^3 = 0 \} \inn
  \fP_{(1,1,1)}=\fP^2. \\ E_2 & = & \{y_0^4+y_1^4+y_2^2 = 0 \} \inn
  \fP_{(1,1,2)}. \\ E_3 & = & \{y_0^6+y_1^3+y_2^2=0\} \inn \fP_{(1,2,3)}.
\end{eqnarray*}
Both $E_1$ and $E_3$ are elliptic curves with modulus $\tau=\grr=e^{2\pi
i/3}$, while $E_2$ has $\tau = i$. 

\subsubsection{K3 surfaces}\label{sK3}
 We let $C_{(w_0,w_1,w_2)}$ denote the
following curve
\begin{equation}\label{11} C_{(w_0,w_1,w_2)} = \{x_0^{\ell} +
  p(x_1,x_2)=0\} \inn \fP_{(w_0,w_1,w_2)}
\end{equation}
of degree $\ell\cdot w_0$. Then applying the map $\Phi$ of Lemma \ref{l3.1}
we get a rational map of $C_{(w_0,w_1,w_2)} \times E_i$ (where $i=1$ or 3
if $\ell=3$ or 6, respectively, 
and $i=2$ for $\ell =4$) onto a hypersurface $X\inn
\fP_{(v_0w_1,v_0w_2,w_0v_1,w_0v_2)}$ of degree $d$, where $(v_0,v_1,v_2)=
(1,1,1),\ (1,1,2)$ and $(1,2,3)$ for the three elliptic curves. We list
some of the possible hypersurfaces, all of which are of Fermat type, 
that one gets in this manner. Assuming that
$d=\sum k_i$, where $(k_1,k_2,k_3,k_4)=(v_0w_1,v_0w_2,w_0v_1,w_0v_2)$, we
necessarily have that $X$ is a (singular) K3 surface, and $\tilde{X}$ is
its minimal desingularisation. The singular fibers we list are those
occuring if the polynomial $p$ in the definition of the curve
$C_{(w_0,w_1,w_2)}$ is of the form $p(x_1,x_2)=x_1^{\ell\cdot w_0/w_1}+
x_2^{\ell\cdot w_0/w_2}$. In the last case the polynomials are not of
Fermat type. 
\[ \begin{array}{|c|c|c|c|c|c|c|} \hline 
\# & (w_0,w_1,w_2) & (v_0,v_1,v_2) & \ell &
  (k_1,k_2,k_3,k_4) & d & \hbox{singular fibers} \\ \hline \hline
1 &  (2,1,1) & (1,1,1) & 3 & (1,1,2,2) &  6 & 6\times IV \\ \hline
2&          & (1,1,2) & 4 & (1,1,2,4) &  8 & 8\times III \\ \hline
3&          & (1,2,3) & 6 & (1,1,4,6) & 12 & 12\times II \\ \hline
4&  (3,1,2) & (1,1,2) & 4 & (1,2,3,6) & 12 & 6\times III, 1\times I_0^* \\ \hline
5&          & (1,2,3) & 6 & (1,2,6,9) & 18 & 9\times II, 1\times I_0^* \\ \hline
6&  (4,1,3) & (1,1,1) & 3 & (1,3,4,4) & 12 & 4\times IV, 1\times IV^* \\ \hline
7&          & (1,2,3) & 6 & (1,3,8,12)& 24 & 8\times II, 1\times IV^* \\ \hline
8&  (5,1,4) & (1,1,2) & 4 & (1,4,5,10)& 20 & 5\times III, 1\times III^* \\ \hline
9&  (7,1,6) & (1,2,3) & 6 & (1,6,14,21)&42 & 7\times II, 1\times II^* \\ \hline
10&  (5,2,3) & (1,2,3) & 6 & (2,3,10,15)&30 & 5\times II, 1\times IV^*,
  1\times I_0^*\\ \hline
11& (11,5,6) & (1,2,3) & 6 & (5,6,22,33) & 66 & 2\times II,\ 2\times II^*
  \\ \hline
\end{array}\]
In a sequel to this paper we will describe in detail how one gets the
singular fibers listed. Let us just make a few remarks about the example
11. For the weights in this case, a Fermat hypersurface is
not possible. We consider instead the following polynomial:
\[
  \{z_0^{12}z_1+z_1^{11}+z_2^3+z_3^2=0\} \inn \fP_{(5,6,22,33)}.\]
We see without difficulty that this is the image under the twist map
\begin{eqnarray*} \fP_{(11,5,6)}\times \fP_{(1,2,3)}  & \lra & 
 \fP_{(5,6,22,33)}\\((x_0:x_1:x_2),(y_0:y_1:y_2)) & \mapsto & 
 (y_0^{5/11}x_1: y_0^{6/11}x_2:x_0^2y_1:x_0^3y_2 ) 
\end{eqnarray*}
of the product $ \{ x_0^6+x_1^{12}x_2 +x_2^{11}
  =0\} \times \{y_0^6+y_1^3+y_2^2=0\}$.
First we note that since $z_1$ only occurs in the
mixed monomial to first power, that the polynomial $p(z_0:z_1:z_2:z_3):=
z_0^{12}z_1+z_1^{11}+z_2^3+z_3^2$ is transversal. It follows that, since
the degree of $X:=\{p(z)=0\}$ is 66, which is also equal to the sum of the
weights, that $X$ is a K3 surface. Its proper transform, after resolving
the singularities, is a smooth K3 surface, which we for convenience also
denote by $X$. From the fact that it is the
image of the twist map above, it follows immediately that 
$X$ has a constant-modulus elliptic fibration.

\subsubsection{An exotic surface}\label{exoticsurface}
 We already mentioned that the condition
$d=\sum k_i$ is sufficient for $X$ to be K3. It is, however, not
necessary, which we show by example. We are indepted to Igor Dolgachev for
explaining this example to us. 
Consider the weighted hypersurface
\[
x^2+y^3 +z^{11} + w^{66} =0\]
in ${\bf P}_{(1,6,22,33)}$. Here $d=66$ while $\sum k_i=62$, hence the
sufficient condition above is not met, and the hypersurface does not look
like a K3. By general principles of weighted hypersurfaces, the geometric 
genus is given by the number of 
monomials of degree $66-62$, of which there is only
one, namely, $w^4$. It follows that a canonical divisor is given by the
locus $w=0$, which is a curve isomorphic to
\begin{equation}\label{12} C=\{z^{11}+y^3+x^2 = 0\} \inn \fP_{(6,22,33)}.
\end{equation} 
The resolution of such a weighted hypersurface is well-known, see
\cite{OW}. One considers
continued fractions of rational numbers derived from the equation of the
hypersurface as follows. Let $(k_1,k_2,k_3,k_4)$ be the weights and for
each $i=1,\ldots,4$, $m_i$ the corresponding degrees of the monomial
$x_i^{m_i}$ occuring in the equation of the hypersurface. Then, to each
singular point, with local stabilizer $\fZ_{\alpha_i}$, the rational number
of which one considers the continued fraction is ${\alpha_i\over \beta_i}$,
where $\beta_i$ is a solution of $k_i\cdot \beta_i = 1mod(\alpha_i)$. In
our case, the $\alpha_i$'s are just the $m_i$'s, and the three fractions we
need to expand are ${2\over 1},\ {3\over 1},\ {11\over 2}$. The first two
give rise to a single $-2$ and $-3$ curve, respectively, while the latter
gives a $-6$ curve meeting a $-2$ curve. Together with the curve $C$ given
here by $x^2+y^3+z^{11}=0$ we have the following resolution of the three
singular points on the K3 surface:
\[ \epsfbox{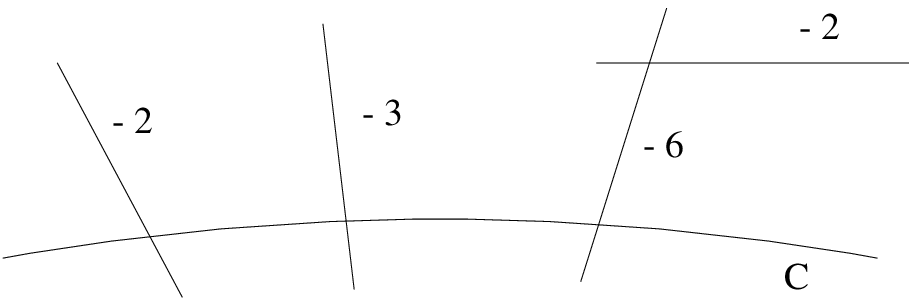} \]
The self-intersection number of the central curve $C$ (the proper transform
of the hyperplane section $w=0$ above) can be determined as follows. Letting
$(\alpha_i,\beta_i)$ denote the pairs giving the continued fractions as
above, the formula is (cf \cite{OW}, 3.6.1, which do the case of an isolated
singularity, but the proof is the same; note that the formula in {\it
  loc.~cit.} should have a minus sign in front of the first term): 
the self-intersection number is
$-b$, where $b$ is given by 
\[ b=-{1\over k_2k_3k_4} + \sum {\beta_i\over \alpha_i}.\]
This yields in our case 
\[ b = -{1\over 66} + {1\over 2} + {1\over 3} + {2\over 11} =1.\]
Hence, the central curve $C$ is {\it exceptional of the first kind} and can
be blown down. The $-2$ curve becomes a $-1$ curve and can be blown down,
then the $-3$ curve of the original curve is now a $-1$ curve and can be
blown down. The result is depicted below. 

\[ \epsfbox{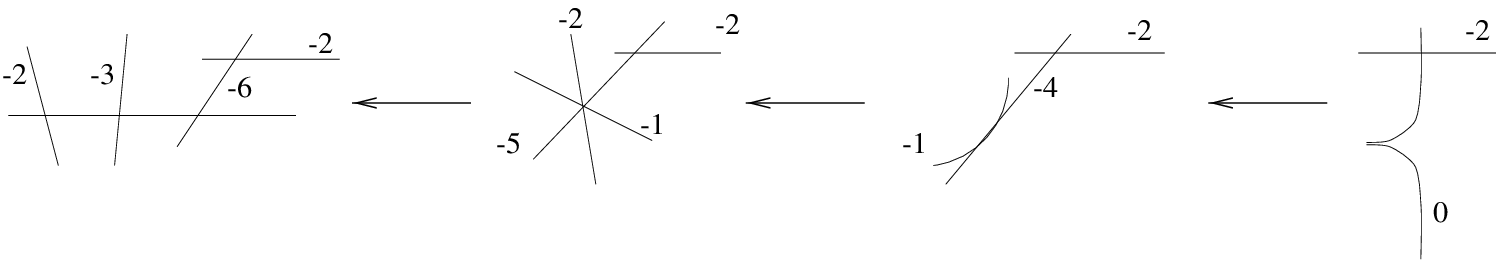} \]
This leaves two curves $F$ and $D$, where
$F$ is a cuspidal genus 1 curve and $D$ is a $-2$ curve intersecting it
transversally. The linear system $|F|$ is an elliptic pencil and $D$ its
section. 
Letting $S$ denote this smooth surface, it is clearly an elliptic K3
surface. It admits an action of $\fZ/66\fZ$ which lifts the action
$(x,y,z,\zeta_{66}w)$ on the weighted hypersurface. 

We can describe this fibration as follows. Consider the subspace $\{z=w=0\}
\cong \fP_{(22,33)}\cong \fP^1
 \inn \fP_{(1,6,22,33)}$. The set of hypersurfaces of
minimal degree containing this $\fP^1$ are given by the equations $D_{\gl}
= \{ z-\gl w^6 =0\}$ (for $\gl\neq 0$). The intersection with the
hypersurface $X$ is then 
\begin{eqnarray*}
 X_{\gl} := D_{\gl}\cap X & = & \{ x^2+y^3+(1+\gl^{11})w^{66} =0\} \cong
\fP_{(1,22,33)}[66] \\ & \cong &  \{ x^2+y^3+(1+\gl^{11})(w')^6 =0\} = 
\fP_{(1,2,3)}[6],
\end{eqnarray*}
which is the elliptic curve of degree 6 in $\fP_{(1,2,3)}$, which is
isomorphic to the Fermat curve, which, as we pointed out above, is the
elliptic curve $E_{\gr}$ with modulus $\gr = e^{2\pi i \over 3}$. The
proper transforms of the 
$X_{\gl}$ are the fibers of our fibration. We need also
the action of $\fZ/66\fZ$ on the holomorphic two-form on the K3
surface. The action on the fibers is:
\[ w \mapsto \zeta_{66} w \Ra w' \mapsto \zeta_{66}^{11} w' =
\zeta_6 w',\]
which is the usual action of $\fZ/6\fZ$ on the sextic in
$\fP_{(1,2,3)}$. The action on the base of the fibration is the action on
the parameter $\gl$ (on some Zariski open set) above, 
\[ \gl w^6 \mapsto \zeta_{11} \gl w^6,\ i.e.,\ z\mapsto \zeta_{11} z,\
\gl\mapsto \zeta_{11} \gl,\]
and the action on the base is just multiplication by $\zeta_{11}$, which
permutes the 11 singular fibers at the 11th roots of $-1$ and fixes the
singular fiber at $\infty$. Since in terms of local coordinates $(z,t)$,
where $z$ is a fiber coordinate, $t$ a coordinate on the base, the
holomorphic two-form is given by $dz\wedge dt$, it follows that 
\[ dz \wedge dt = dw \wedge d\gl \mapsto \zeta_6 dw \wedge \zeta_{11} d\gl
= \zeta_{66} dz\wedge dt,\]
and the action is by multiplication by $\zeta_{66}$, i.e., if $h$ denotes
the generator of $\fZ/66\fZ$ and $\go$ denotes the holomorphic two-form,
then 
\begin{equation}\label{e66} h^*(\go) = \zeta_{66} \go.
\end{equation}

This surface is known
from the work of Kondo. He shows in \cite{Kondo2} that there is in fact a
{\it unique} K3 surface which admits $\fZ/66\fZ$ as an automorphism
group which preserves the Picard group of the K3 surface, or equivalently,
acts non-trivially on the holomorphic two-form (see section 3.3.2 below). 
Kondo describes the same surface as an elliptic surface 
with 12 fibers of type \II\ 
at $t=0$ and at the $11^{th}$
roots of unity. The affine equation is 
\[y^2=x^3+t\prod_1^{11}(t-\zeta_{11}^i),\]
and the automorphism is given by 
\[ g: (x,y,t) \mapsto (\zeta_{66}^2 x, \zeta_{66}^3 y, \zeta_{66}^6 t).\]
He shows that the Picard lattice of this surface is just a hyperbolic
plane, while the transcendental lattice is of the form $U\oplus U \oplus E_8
\oplus E_8$, where $E_8$ denotes the negative definite copy of the root
lattice of the exceptional group of type $E_8$. For the elliptic surface in
Weierstra\ss\ form, the holomorphic two-form is given by ${dx\over y}
\wedge dt$, and it follows that $g^*(\go)=\zeta_{66}^5 \go$, which is a
primitive $66$th root of unity, showing that the action of the group of
order 66 acts trivially on the Picard group. Comparing with (\ref{e66}), we
see that $h=g^{55}$. 

\subsection{Calabi-Yau threefolds}
In this section we apply the twist map to construct Calabi-Yau threefolds
with K3-fibrations, which are quotients of either products $S\times E$,
where $S$ is a surface and $E$ is an elliptic curve, or products 
$C\times K$, where $K$
is now a K3 surface. To apply the latter, we need some information on the
automorphisms of K3 surfaces. This information can be found in
\cite{Kondo2} and \cite{Kondo1}. We first recall this, then pass to the
construction. 

\subsubsection{Elliptic fibrations}\label{selliptic}
 We do not attempt to give a
classification, but instead describe a few examples. Lists of such elliptic
weighted hypersurface Calabi-Yau threefolds have been compiled, and for
many of these we can realize them as quotients. As before we use the three
elliptic curves of section 3.2.1 but instead of the curve (\ref{11}) now
the variety $V_1$ is the surface 

\begin{table}\[ \begin{array} {|l|l| c | l | c|} \hline (w_0,w_1,w_2,w_3,)
    & (v_0,v_1,v_2) & \ell & (k_1,k_2,k_3,k_4,k_5) & d \\ \hline \hline 
    (3,1,1,1) & (1,1,1) & 3 & (1,1,1,3,3) & 9 \\ 
              & (1,1,2) & 4 & (1,1,1,3,6) & 12 \\ 
              & (1,2,3) & 6 & (1,1,1,6,9) & 18 \\ \hline 
    (4,1,1,2) & (1,1,1) & 3 & ( 1,1,2,4,4) & 12 \\
              & (1,1,2) & 4 & (1,1,2,4,8) & 16 \\ 
              & (1,2,3) & 6 & (1,1,2,8,12) & 24 \\ \hline
    (5,1,1,3) & (1,1,1) & 3 & (1,1,3,5,5) & 15 \\
              & (1,2,3) & 6 & (1,1,3,10,15) & 30 \\ \hline
    (5,1,1,2) & (1,1,2) & 4 & (1,2,2,5,10) & 20 \\ 
              & (1,2,3) & 6 & (1.2.2.10.15) & 30 \\ \hline 
    (6,1,1,4) & (1,1,2) & 4 & (1,1,4,6,12) & 24 \\
              & (1,2,3) & 6 & (1,1,4,12,18) & 36 \\ \hline
    (6,1,2,3) & (1,1,1) & 3 & (1,2,3,6,6) & 18 \\
              & (1,1,2) & 4 & (1,2,3,6,12) & 24 \\
              & (1,2,3) & 6 & (1,2,3,12,18) & 36 \\ \hline
    (7,1,2,4) & (1,1,2) & 4 & (1,2,4,7,14) & 28 \\ \hline
    (7,1,3,3) & (1,1,1) & 3 & (1,3,3,7,7) & 21 \\ 
              & (1,2,3) & 6 & (1,3,3,14,21) & 42 \\ \hline
    (7,2,2,3) & (1,2,3) & 6 & (2,2,3,14,21) & 42 \\ \hline
    (8,1,1,6) & (1,1,1) & 3 & (1,1,6,8,8) & 24 \\ 
              & (1,2,3) & 6 & (1,1,6,16,24) & 48 \\ \hline
    (8,1,3,4) & (1,2,3) & 6 & (1,3,4,16,24) & 48 \\ \hline
    (8,2,3,3) & (1,2,3) & 6 & (2,3,3,16,24) & 48 \\ \hline
    (9,1,2,6) & (1,1,2) & 4 & (1,2,6,9,18) & 36 \\ 
              & (1,2,3) & 6 & (1,2,6,18,27) & 54 \\ \hline
    (9,1,4,4) & (1,1,2) & 4 & (1,4,4,9,18) & 36 \\ \hline
    (9,2,3,4) & (1,1,2) & 4 & (2,3,4,9,18) & 36 \\ \hline
    (10,1,1,8)& (1,1,2) & 4 & (1,1,8,10,20) & 40 \\ \hline
    (10,2,3,5) & (1,1,1) & 3 & (2,3,5,10,10) & 30 \\
               & (1,2,3) & 6 & (2,3,5,20,30) & 60 \\ \hline
    (10,1,3,6) & (1,1,1) & 3 & (1,3,6,10,10) & 30  \\ 
               & (1,2,3) & 6 & (1,3,6,20,30) & 60 \\ \hline
    (10,3,3,4) & (1,2,3) & 6 & (3,3,4,20,30) & 60 \\ \hline
    (12,1,2,9) & (1,1,1) & 3 & (1,2,9,12,12) & 36 \\
               & (1,2,3) & 6 & (1,2,9,24,36) & 72 \\ \hline
    (13,1,6,6) & (1,2,3) & 6 & (1,6,6,26,39) & 78 \\ \hline
    (14,1,1,12) & (1,2,3) & 6 & (1,1,12,28,42) & 84 \\ \hline

\end{array}\]\caption{List of elliptic Calabi-Yau threefold
    weighted hypersurfaces of Fermat type and with constant fiber modulus}
          \end{table}

\begin{equation}\label{surface} S_{(w_0,w_1,w_2,w_3)} =\{ x_0^{\ell} +
  p(x_1,x_2,x_3) = 0 \}
\end{equation}  
of degree $\ell w_0$. Applying the map $\Phi$ of Lemma \ref{l3.1} we get a
rational map of $S_{(w_0,\ldots, w_3)} \times E_i$ (where $i=1$ or 3
if $\ell=3$ or 6, respectively, and $i=2$ for $\ell =4$) onto a threefold of
degree $d=\ell w_0$. again we list some of the possiblities which a brief
manual search comes up with. All these cases have the following properties:
they are of Fermat type (quite restrictive) and they have a fibration of
{\it constant modulus} (it is not clear how restrictive this is). 
This list is given in Table 1.

As an explicit example consider the case $(2k,1,1,2(k-1)),\ (1,2,3)$ for
which the Calabi-Yau is the hypersurface of degree $12k$ in
$\fP_{(1,1,2(k-1),4k,6k)}$. The elliptic fibration is onto
$\fP_{(1,1,2(k-1))}$, whose desingularisation is the Hirzebruch surface
$\fF_{2(k-1)}$. By Lemma \ref{l3.9}, the discriminant is the total
transform of $\gD=\{p(x_1,x_2,x_3)=0\} \inn \fP_{(1,1,2(k-1))}$. The
projection $\fP_{(1,1,2(k-1))} \lra \fP_{(1,1)}$ lifts to the projection
$\fF_{2(k-1)} \lra \fP^1$. Let $C_{\infty}$ be the section of negative
self-intersection $-2(k-1)$ (the exceptional curve of the resolution
$\fF_{2(k-1)} \lra \fP_{(1,1,2(k-1))}$), and let $C_0$ denote the class of
a positive section, $F$ the class of a fiber. The curve $\gD$ is reducible;
its total transform on $\fF_{2(k-1)}$ consists of an irreducible curve $C$
of class $C\sim aC_0+bF$ ($a,b>0$) and a multiple of $C_{\infty}$, $\gD
\sim a C_0 + bF + \nu C_{\infty}$. Assume that the exceptional point
$(0,0,1)$ of $\fP_{(1,1,2(k-1))}$ is not contained in $\gD$ (for example
$p$ of Fermat type). Then $C$ and $C_{\infty}$ are disjoint. Hence given
the explicit form of $p$ one can easily determine exactly the degeneracy
locus of the smooth elliptic Calabi-Yau.

In the cases $k=2,3,4,7$ and $\ell=6$ or $k=5$ and $\ell =4$
 one can take $p$ to be of Fermat type. In these
cases, the class of $C$ can be calculated as follows: it is a curve in
$\fP_{(1,1,2(k-1))}$ of degree $12k$, 
which maps $k$ to one onto a $\fP^1$, totally branched at $12k$ points,
hence the Euler number (which is the first Chern number) is equal to
$k\cdot(2-2k)+2k=2k(2-k)$. Then applying adjunction on the resolution
$\fF_{2(k-1)}$, one gets $a=k $ and $b=-2k(k-1) $. The fiber over
$C_{\infty}$ is of type IV, I$_0^*$, IV$^*$ and II$^*$ in the $\ell=6$
cases and III$^*$ in the $k=5,\ \ell=4$ case, and this determines $\nu$ to
be 3, 6, 8, 10 and 9 in the respective cases. These have been studied in
more detail in \cite{MV}. 

\begin{lemma}\label{l5} Let $X$ be a Calabi-Yau threefold with both an
  elliptic and a K3 fibration, and suppose these are compatible (i.e., that
  the elliptic fibration of $X$ restricts to an elliptic fibration of the
  smooth fibers), and assume moreover that the degenerate fibers of the K3
  fibration are irreducible. Then the base of the elliptic fibration is a
  rational ruled surface $\bf F_n$. 
\end{lemma}
{\bf Proof:} From the first assumption, we clearly have on the base of the
elliptic fibration, for each K3 fiber, a $\fP^1$, the base of the elliptic
fibration of that K3 fiber. Thus the surface fibers over the original
$\fP^1$, and is thus the blow up of one of the Hirzebruch surfaces. Now we
use the assumption that each of the degenerate fibers is irreducible; being
so, each of these can fiber over at most a single $\fP^1$, and it follows
that {\it all} fibers of the base surface are just $\fP^1$, in other words
that it is a rational ruled surface, as claimed. \ende

In general the base of such a fibration can be some blow-up of one of the
$\fF_n$. 
\begin{table}
\[ \begin{array}{|l|l||c|c|c|} \hline
(w_0,w_1,w_2,w_3)  & (k_1,k_2,k_3,k_4,k_5) & d & \chi &
h^{1,1} \\ \hline \hline
(581,41,42,498) & (41,42,498,1162,1743) & 3486 & 960 & 491 \\ \hline
(498,36,41,421) & (36,41,421,996,1494) & 2988 & 960 & 491 \\ \hline
(539,36,41,462) & (36,41,462,1078,1617) & 3234 & 900 & 462 \\ \hline
(469,31,42,396) & (31,42,396,938,1407) & 2814 & 900 & 462 \\ \hline
(463,31,41,391) & (31,41,391,926,1389) & 2778 & 900 & 462 \\ \hline
(433,31,36,366) & (31,36,366,866,1299) & 2598 & 840 & 433 \\ \hline
(483,28,41,414) & (28,41,414,966,1449) & 2898 & 804 & 416 \\ \hline
(414,24,41,349) & (24,41,349,828,1242) & 2484 & 804 & 416 \\ \hline
(385,28,31,326) & (28,31,326,770,1155) & 2310 & 744 & 387 \\ \hline
(434,21,41,372) & (21,41,372,868,1302) & 2604 & 720 & 377 \\ \hline
(372,18,41,313) & (18,41,313,744,1116) & 2232 & 720 & 377 \\ \hline
\end{array}\]
\caption{\label{table2}
  Cases of elliptic fibrations with constant modulus (fiber the
  sextic elliptic curve in $\fP_{(1,2,3)}$), for which the Euler number
  $\chi$ is large, positive}
\end{table}

We now turn to some examples of elliptic fibrations which are not K3
fibrations. The examples given up to now were all in the geography region
where the Euler number is negative. We now give some examples where the
Euler number is positive, which are listed in Table \ref{table2}. 
We will discuss one example in more detail.
These examples are all images of the twist map 
\[ \fP_{(w_0,w_1,w_2,w_3)} \times \fP_{(1,2,3)} \lra
\fP_{(w_1,w_2,w_3,2\cdot w_0,3\cdot w_0)}.\]
Hence they all have constant elliptic fibrations with fiber the sextic in
$\fP_{(1,2,3)}$. We now consider the first example in some more detail.

The variety $V_1$ is a weighted hypersurface in $\fP_{(581,41,42,498)}$ of
degree 3486, given by the equation 
$$x_0^6+p(x_1,x_2,x_3) =  x_0^{6}+x_1^{84}x_2 + x_2^{83} +
x_3^{7} = 0$$ and is a $6-1$ cover of the weighted $\fP_{(41,42,498)} \cong
\fP_{(41,7,83)}$, branched over the locus $p=0$. Then the base of the
elliptic fibration will be the resolution of $\fP_{(41,7,83)}$. This is
standard toric geometry. Take the lattice spanned by the three vectors
\[ v_0 = {1\over 41} \left( \begin{array}{c} -1 \\ -1  
  \end{array} \right), \quad v_1={ 1\over 7} \left(\begin{array}{c} 1 \\ 0
  \end{array} \right), \quad v_2={1\over 83}\left(\begin{array}{c} 0 \\ 1 
  \end{array} \right), \]
and in each of the cones spanned by two of the $v_i$ one must determine all
lattice points. Then taking the convex hull of these, one gets a simplicial
decomposition of the cones into subcones all of which have unit area. Then
one counts the number of new vertices, and each corresponds to an
exceptional divisor. There are three singular points, and our curve $p=0$
passes through only $(1,0,0)$, so this is the only relevant singularity. 
It is easy to see that we are looking for integral solutions
$(\ga,\gb,\gg)$ of the
inequalities
\[ \ga(w_1+w_2)-w_0(\gb+\gg) + w_0 \leq 0,\quad w_0\gb>w_1\ga,\quad
w_0\gg>w_2\ga,\]
where $(w_0,w_1,w_2)=(41,7,83)$, 
of which there are the following 20, 
in other words, there are 20 exceptional curves
resolving the point: 
\begin{eqnarray*} (\ga,\gb,\gg) & = & (11, 2, 23), (16, 3, 33), (17, 3,
35), (21, 4, 43), (22, 4, 45), (23, 4, 47), (26, 5, 53), \\ & & 
 (27, 5, 55), (28, 5, 57), (29, 5, 59), (31, 6, 63), (32, 6, 65, )(33, 6,
 67),  (34, 6, 69), \\ & & 
(35, 6, 71), (36, 7, 73), (37, 7, 75), (38, 7, 77), (39, 7, 79), (40, 7, 81).
\end{eqnarray*}
It follows that along the proper transform of $p=0$,
we have singularities of type II, while along the 20 exceptional curves we
have singularities of type II$^*$. The morphism onto the base is not yet
flat; to achieve this one could take Miranda's small resolutions, which
would amount to blowing up the intersection points of the 20 exceptional
curves.

\subsubsection{Automorphisms of K3 surfaces} 
This section is prepatory and recalls some facts about the automorphism of
K3 surfaces, which will be applied in the same way the automorphisms of
orders 2, 3, 4 and 6 were for elliptic curves.
Let $X$ be a K3 surface; the integral homology ${\rm H}^2(X,\fZ)$ forms a lattice
in ${\rm H}^2(X,\fR)\cong \fR^{22}$ which is a copy of $U\oplus U \oplus U \oplus
E_8 \oplus E_8$, where $E_8$ denotes here the negative definite even
unimodular lattice in $\fR^8$ and $U$ is the hyperbolic lattice. We
consider in this note {\it polarized} K3 surfaces, letting $\omega \in
{\rm H}^{1,1}(X)$ denote the K\"ahler form; the primitive cohomology
${\rm H}^2(X,\fC)_0$ is then the orthocomplement in ${\rm H}^2(X,\fC)$ to
$<\omega>$. By an {\it automorphism} of $X$ we will mean in this paper a
polarization-preserving automorphism $\sigma:X\lra X$; we let also
$\sigma_*$ denote the induced map on cohomology $\sigma_*\in
\Aut({\rm H}^2(X,\fZ))$. The map 
\begin{eqnarray*} *: \Aut(X) & \lra & \Aut({\rm H}^2(X,\fZ)) \\ \sigma &
  \mapsto & \sigma_*
\end{eqnarray*}
is injective (cf. \cite{Kondo1}, section 4) so it is convenient to consider
this as a representation of $\Aut(X)$. 
Let $S_X:={\rm H}^{1,1}(X) \cap {\rm H}^2(X,\fZ)$ be the {\it Picard lattice} of $X$,
and let $T_X$ be its orthocomplement in ${\rm H}^2(X,\fZ)$; this is the {\it
  transcendental lattice} of $X$. The direct sum $S_X\oplus T_X$ is of
finite index in ${\rm H}^2(X,\fZ)$, and the map $*$ actually maps into
$\Aut(S_X\oplus T_X) \cong O(S_X)\times O(T_X)$ (cf.\ \cite{Kondo2},
Proposition 1.1), so we can faithfully consider the action of $\Aut(X)$ on
the Picard and transcendental lattices. Note that if $T_X$ (and hence
$S_X$) is {\it unimodular}, then 
\[  S_X\oplus T_X={\rm H}^2(X,\fZ)=U^{\oplus 3}\oplus E_8^{\oplus 2},\]
so both of $S_X$ and $T_X$ are themselves direct sums with summands $U$ and
$E_8$. 

Every element $g\in \Aut(X)$ acts on this
decomposition, and letting $\alpha_g\in \fC^*$ denote the factor
such that $g^*(\Omega)=\alpha_g \Omega$, where $\Omega$ denotes the
holomorphic two form, the map $g\mapsto \alpha_g$ 
gives rise to an exact sequence 
\[ 1 \lra G_X \lra \Aut(X) \stackrel{\alpha}{\lra} H_X \lra 1.\]
Elements of $G_X$ are called {\it symplectic}, as they preserve the
symplectic form $\Omega$. If $X$ is algebraic, then $H_X=\fZ_k$ for some
$k$. The following facts were proved by Nikulin in \cite{vvn}:
\begin{itemize}\item[1.] $G_X=\Ker(\Aut(X) \lra \Aut(T_X))$ and
  $H_X=\Ker(\Aut(X)\lra \Aut(S_X))$. \item[2.] The representation
  $\fZ_k\lra \Aut(T_X)$ is a direct sum of a number 
  $K$ of irreducible representations,
  each of rank $\phi(k)$ ($\phi$ denotes the Euler phi-function), 
  so $\sum^K_1\phi(k)=K \phi(k) =rank(T_X)$. \item[3.]
  The representation $G_X\lra \Aut(S_X)$ was described by Nikulin for
  abelian $G_X$; the possible subgroups occuring are  
  $$({\bf Z}_2)^m, \ 0\leq m\leq 4;\ \ {\bf Z}_4;\ \ {\bf Z}_2\times
  {\bf Z}_4;\ \ ({\bf Z}_4)^2;\ \ {\bf Z}_8;\ \ {\bf Z}_3;$$
   $$ ({\bf
    Z}_3)^2;\ \ {\bf Z}_5;\ \ {\bf Z}_7; \ \ {\bf Z}_6;\ \ {\bf Z_2}\times
  {\bf Z}_6.$$
\end{itemize}
From the fact that $\phi(k)\leq 22$ one deduces $k\leq 66$, and $k=66$ does
occur. This was considered by Kondo in \cite{Kondo2} and by Vorontsov in 
\cite{Vor}, who proved the following results:
\begin{itemize}\item[1.] If $T_X$ is unimodular, then $k$ is a divisor of
  66, 44, 42, 36, 28 or 12. \item[2.] Suppose $\phi(k)=rank(T_X)$.  
  Then $k$=66, 44, 42, 36, 28 or 12,
  and in these cases there exists a unique (up to isomorphism) K3 surface
  with given $k$. \item[3.] If $T_X$ is not unimodular, then $k=2^r\ (1\leq
  r \leq 4),\ 3^r\ (1\leq r \leq 3),\ 5^r\ (r=1, 2),\ 7, 11, 13, 17$ or 19. 
\end{itemize}
In the case that $T_X$ is unimodular, it is a direct sum of factors $U$ and
$E_8$. For the cases above, the lattices are listed in the following
table.
\[\begin{array}{c | c c} k & S_X & T_X \\ \hline
66 & U & U\oplus U \oplus E_8 \oplus E_8 \\
44 & U &  U\oplus U \oplus E_8 \oplus E_8 \\
42 & U\oplus E_8 &  U\oplus U \oplus E_8  \\
36 & U\oplus E_8 &  U\oplus U \oplus E_8  \\
28 & U\oplus E_8 &  U\oplus U \oplus E_8  \\
12 & U\oplus E_8 \oplus E_8 & U \oplus U 
\end{array}
\]

\subsubsection{Calabi-Yau threefolds with K3 fibrations}
Again we just list some examples. As above we consider three K3's of the
above list. 
\begin{eqnarray*} K_1 & = & \{y_0^6+y_1^6+y_2^3+y_3^3 =0\}\inn
  \fP_{(1,1,2,2)} \\
                  K_2 & = & \{y_0^{12}+y_1^6+y_2^4+y_3^2=0\} \inn
                  \fP_{(1,2,3,6)} \\
                  K_3 & = & \{y_0^{42}+y_1^7+y_2^3+y_3^2 =0 \} \inn
                  \fP_{(1,6,14,21)}.
                \end{eqnarray*}
All three K3's are elliptic fibrations;                
the elliptic fibers are $E_1,\ E_2$ and $E_3$, respectively, and $\ell = 6,
12, 42$. We again consider the curve (\ref{11}), and the weighted
hypersurfaces are then birational to quotients of $C_{(w_0,w_1,w_2)}\times
K_i$ by $\fZ/\ell\fZ$. We list some examples found after a brief manual
search in Table 3. Of course most of these also occur in Table 1.

\begin{table}
\[ \begin{array}{|l|l| c | l | c|c|} \hline (w_0,w_1,w_2)
    & (v_0,v_1,v_2,v_2) & \ell & (k_1,k_2,k_3,k_4,k_5) & d & \chi
    \\ \hline \hline
 (2,1,1) & (1,1,2,2) & 6 & (1,1,2,4,4) & 12 & -192 \\
         & (1,2,3,6) & 12 & (1,1,4,6,12) & 24 & -312 \\
         & (1,6,14,21) & 42 & (1,1,12,28,42) & 84 & -960 \\ \hline
 (3,1,2) & (1,1,2,2) & 6 & (1,2,3,6,6) & 18 & -144 \\
         & (1,2,3,6) & 12 & (1,2,6,9,18) & 36 & -228 \\
         & (1,6,14,21) & 42 & (1,2,18,42,63) & 126 & -720 \\ \hline
 (4,1,3) & (1,1,2,2) & 6 & (1,3,4,8,8) & 24 & -120 \\ 
         & (1,2,3,6) & 12 & (1,3,8,12,24) & 48 & -192 \\
         & (1,6,14,21) & 42 & (1,3,24,56,84) & 168 & -624 \\ \hline
 (5,1,4) & (1,2,3,6) & 12 & (1,4,10,15,30) & 60 & -168 \\ \hline
 (7,1,6) & (1,2,3,6) & 12 & (1,6,14,21,42) & 84 & -132 \\
         & (1,6,14,21) & 42 & (1,6,42,98,147) & 294 & -480 \\ \hline
 (5,2,3) & (1,1,2,2) & 6 & (2,3,5,10,10) & 30 & -72 \\
         & (1,2,3,6) & 12 & (2,3,10,15,30) & 60 & -108 \\
         & (1,6,14,21) & 42 & (2,3,30,70,105) & 210 & -384 \\ \hline
       \end{array}
        \]\caption{K3-fibered Calabi-Yau weighted hypersurfaces which are
 also elliptic fibered, have constant modulus and are of Fermat type}
      \end{table}

We now show by example that {\it a Calabi-Yau threefold can have two
  different K3-fibrations with constant modulus}. 
\begin{lemma}\label{theorem2} An appropriately choosen
  weighted hypersurface of degree 12 in the weighted
  projective space $\fP_{(2,2,1,1,6)}$ 
  has two K3-fibrations with constant modulus; the two
  different fibers are $\fP_{(1,1,4,6)}[12]$ and $\fP_{(1,1,1,3)}[6]$.
\end{lemma}
{\bf Proof:} Both of these fibrations can be constructed with the twist
map. For the first, we let $\boldmu_3$ act on the product
$\fP_{(2,1,1)}[6]\times \fP_{(4,1,1,6)}[12]$; the twist map is onto a
hypersurface of degree 24 in $\fP_{(4,4,2,2,12)}$, which is the same thing
(as all weights are divisible by 2) as a hypersurface of degree 12 in
$\fP_{(2,2,1,1,6)}$. By Corollary \ref{cfiber}, we get a constant modulus
fibration with fiber $\fP_{(4,1,1,6)}[12]$. For the second, we let
$\boldmu_6$ act on the product $\fP_{(2,1,1)}[12]\times
\fP_{(1,1,1,3)}[6]$. The image under the twist map is a hypersurface of
degree 12 in $\fP_{(2,2,1,1,6)}$, and by choosing equations of the K3
surface and auxiliary curves appropriately, these two hypersurfaces
coincide. In particular it is true of the Fermat hypersurface. \ende 

\subsubsection{A curious example}

Consider the ten K3 surfaces of Fermat type in the table in section
\ref{sK3}. For each, of weights $(k_1,k_2,k_3,k_4)$ and degree $d=\sum
k_i$, we can form the image under the twist map with $\fP_{(2,1,1)}$, 
\[ \fP_{(2,1,1)}[2d] \times \fP_{(k_1,k_2,k_3,k_4)}[d] \lra
\fP_{(1,1,2k_2,2k_3,2k_4)}[2d].\] For the cases 1, 4, and 9 in that list, the
resulting Calabi-Yau threefolds are the first three entries in Table
3. The analysis described there applies also to these cases, and the K3
fibration over $\fP^1$ has $2d$ singular fibers, which are the affine K3
singularities listed in Table 2 of part II of this paper. 
Let {\bf X} be the type of bad fiber,
$e({\bf X})$ its Euler number. Then we have the formula for the Euler
number of the Calabi-Yau
\[ e(X) = (2-2d) \cdot 24 + 2d \cdot e({\bf X}) = 48 + 2d\cdot (e({\bf X})
- 24).\]
Now since $0<e({\bf X})<24$, we have $-24 < e({\bf X}) -24 <0$ which yields
the inequality for the Euler number of $X$:
\[ -48 d < e(X) - 48 = 2d(e({\bf X})-24) < 0.\]

Observe that $\fZ/d\fZ$ is an automorphism group of the K3 $V_2$,
so by the results described above we have $d\leq 66$. Realizing that the case
with minimal known Euler number $-960$ is realized in this manner by taking
$V_2$ to be the K3 with $d=42$ (this is the third example in Table 3), 
while there are cases $d=44, 66$, one could
imagine constructing a similar example with the $d=66$ case. 

Next observe that if $g$ denotes the generator of $\Aut(V_2)$, fixing $x\in
\gD\inn \fP^1$ in the discriminant of the fibration, formula (\ref{8})
shows the contribution of the singular fiber $X_x$ at $x\in \gD$ is
determined by the action of $g$ on the holomorphic two-form. In each case
of the fibrations described above, for each $x_i\in \gD,\ a_i={1\over d}$,
and we have the equality 
\[ \sum_{i=1}^{2d} a_i = 2d {1\over d} = 2 = c_1(\fP^1),\]
which is the necessary condition described at the beginning of this
paper. For the case $d=66$ this would require 132 singular fibers, each
giving a contribution of ${1\over 66}$. If so, we can calculate what its
Euler number would be. Since it
has 132 singular fibers which are the affine surface $t_1^{11} +t_2^3
+t_4^2 = 0 $ which has Milnor number $\mu = (11-1)(3-1)(2-1)=20$, each
singular fiber has Euler number 4, and our formula for the Euler number
yields:
\[ e(X) = (2-132) \cdot 24 + 132\cdot 4 = -2592.\]
Is it possible to contruct such an
example? 

Consider the image of the twist map 
\[ \fP_{(2,1,1)}[132] \times \fP_{(1,6,22,33)}[66] \lra
\fP_{(1,1,12,44,66)}[132].\]
If we take $V_1,\ V_2$ as Fermat hypersurfaces, then the image is 
the weighted threefold given by the equation 
\[ X = \{ z_1^{132} + z_2^{132} + t_1^{11} +t_2^3 +t_4^2 = 0 \} \inn
\fP_{(1,1,12,44,66)}. \]
Suppose this threefold did allow a fibration with fibers the resolved K3
surfaces discussed in section \ref{exoticsurface}. 
It is easy to see that we have
singular fibers at the 132 points which are the zeroes of $z_1^{132} +
z_2^{132}$, and that the fibers $\{z_1=0\}\cap X$ and $\{z_2=0\}\cap X$ are
both smooth, so that it looks as if we could in fact construct such a
fibration. However, as $132 > 124=\sum k_i$, the threefold $X$ does not
satisfy the sufficient condition to be Calabi-Yau, and our arguments above
(in particular Lemma \ref{l3.4})
do not apply directly. In fact, assuming the necessary condition (\ref{4})
is also sufficient, such a fibration does {\it not} exist,
which can be seen as follows. There is no birational model of $X$ which is
Calabi-Yau: the geometric genus
$h^{3,0}$ is 9, not 1 as would be the case if $X$ were Calabi-Yau. Indeed, 
since the geometric genus is a birational invariant, this is
just the number of sections of $\cO(132-124)=\cO(8)$, and this is the
number of monomials in the first two variables of weight 8, of which there
are 9. However, using the method of the proof of Lemma \ref{l3.4}, the
reader may verify that in fact such a fibration does exist (first the fiber
is the non-Calabi-Yau surface $\{x^2+y^3+z^{11}+w^{66}=0\}\inn
\fP_{(1,6,22,33)}$ of section 3.3.3, which is birationally modified as in
that section over the base of the fibration). Consequntly, condition
(\ref{4}) is not sufficient. 

We note that we can also display 
$X$ as an elliptic fibration, and for these we have the necessary and
sufficient condition (\ref{6}), which utilizing the Weierstra\ss\ form says
that 
\[ -12 K_Y = \gD,\]
where $\gD$ is the discriminant locus (counted with appropriate
multiplicities). The variety $X$ is the image of the twist map
\[ \Phi:\fP_{(11,1,1,12)}[132] \times \fP_{(1,2,3)}[6] \lra
\fP_{(1,1,12,44,66)}[132].\]
The discriminant is the total transform of 
\[ \gS = \{x_1^{132}+x_2^{132}+x_3^{11} =0 \} \inn \fP_{(1,1,12)};\]
projecting onto $\fP_{(1,1)}$ displays the {\it proper transform} of this,
which we denote also by $\gS$, as an 11 to 1 cover, totally
branched at the 132rd roots of $-1$, so it has Euler number
\[ e(\gS) = 11(2-132) + 132 = -1298.\]
We can determine its class in ${\rm H}^2(\fF_{12},\fZ)$ as follows: from the fact
that it is an 11 to 1 cover, it intersects a fiber in 11 points, so 
$\gS \sim 11 C_0 + b F$, where $C_0$ denotes the (class of) a section of
positive self-intersection and $F$ denotes the class of a fiber. On the
other hand, the first Chern class of the resolution $\fF_{12}$ of
$\fP_{(1,1,12)}$ is 
\[ c_1(\fF_{12}) = 2 C_0 - 10 F,\]
and applying adjuction to $\gS$ on $\fF_{12}$,
which is a smooth curve with Euler number we just calculated, one can
determine its class in ${\rm H}^2(\fF_{12},\fZ)$. The result is 
\[ \gS \sim 11 C_0 \sim 11( C_{\infty} +12F),\]
while for the exceptional curve $C_{\infty}$, we have fiber types II again,
so 
\[ \gD = 2\gS + 2 C_{\infty} = 24 C_{\infty} + 264 F,\]
and one sees that the sufficient condition (\ref{4}) is not
satisfied, as $-12K_{\fF_{12}} = 24 C_{\infty} + 168 F$. 

It remains an open problem whether it is possible to construct such a K3
fibration; if possible, it would enlarge the range of possible Euler
numbers of Calabi-Yau threefolds (and is therefore quite unlikely). 

\subsubsection{Birational fibrations}
It is, however, possible to construct K3 fibrations which have the K3
surface of section \ref{exoticsurface} as fiber, which we now explain.
There are examples of weighted hypersurfaces which, after resolution of
singularities, are not fibrations, but still are, in the class of
Calabi-Yau threefolds, birational to such a Calabi-Yau. These examples come
from the exotic surface example of section \ref{exoticsurface}. Indeed,
suppose we want a twist map 
\[ \Phi : \fP_{(w_0,w_1,w_2)}[66\cdot w_0]\times \fP_{(1,6,22,33)}[66] \lra
\fP_{(w_1,w_2,w_0\cdot 6,w_0\cdot 22,w_0\cdot 33)}[66\cdot w_0]\]
to have an image which satisfies the sufficient condition $d = \sum k_i$
for it to be a Calabi-Yau threefold. 
Then writing down what the weights are,
we get an equation $w_1+w_2 +w_0(6+22+33) \stackrel{!}{=}
 66w_0$, which means that we 
require $w_1+w_2 =5w_0$. The easiest solution to this is: 
$w_0=1, (w_1,w_2) = (2,3)$ or 
$(1,4)$. The image projective spaces is then $\fP_{(2,3,6,22,33)}$, and 
there are indeed Calabi-Yaus (of degree 66) in this space, as well as in 
$\fP_{(1,4,6,22,33)}$. They have Euler numbers $-240$ and $-300$, 
respectively. 
They do not posess a fibration {\it a priori}, because the base locus 
of the projection (which one gets by setting the 
first two coordinates =0) is not part of the singular 
locus (recall that in Lemma \ref{l3.4} we required $w_0>1$ to get a
fibration). However, one can blow it up (upon which the surface is temporarily
no longer C-Y) and then in each fiber do the blowing down process we have
described in section \ref{exoticsurface}.
After this is done, we do get a C-Y with fibration by those exotic surfaces.

\section{Applications in Physics}

In this section we breifly describe some of the applications of the twist
map to physical dualities, which was the original source of motivation for
the present investigation. It is intended for the non-expert, and we just
try to explain the physical interpretation of the geometry, without going
into any details.

\subsection{The physical theories and their moduli spaces}

We are interested in {\it superstring theories}. These are theories about
how a string (a smooth image of the circle) moves in Minkowski space; its
vibrations give rise to all elementary particles. More precisely, consider
a string moving in some Minkowski space $M^{1,d}$; it traces out with time
a {\it world sheet}, and this is then an embedded Riemann surface in
$M^{1,d}$. The physical theory one is interested in is a {\it
  superconformal field theory} on that Riemann surface; the adjective
superconformal refers to the symmetry group of the equations of motion. It
is well-known that for this symmetry to hold, the dimension is restricted
to $d=9$ (ten-dimensional Minkowski space), and it is also known that there
are five consistent theories of this type: Types I, IIA, IIB,
Het$_{E_8\times E_8}$, Het$_{SO(32)}$. The first is the theory which
contains open strings and has only $N=1$ supersymmetry, while the others
are theories of closed strings and have $N=2$ supersymmetry (on the world
sheet). That is, one has a consistent supersymmetric 
quantum field theory on the world
sheet of the string. (There is also an issue of space-time supersymmetry,
but we neglect that here). The type II theories have at most abelian gauge
groups (at least in ten dimensions, see \cite{SS2}, Chapter 14 for a
discussion of this issue upon compactification), 
which is why they seemed completely uninteresting for many years,
while the heterotic strings, which have only left-moving supersymmetries,
have the two gauge groups of rank 16 (which arises from the fact that the
consistent dimension for the {\it bosonic}, or non-supersymmetric string,
is 26, and the $16=26-10$ remaining dimensions are compactified to the
maximal torus of a Lie group, the gauge group) listed above. 

The superstring theories are all {\it perturbative} in essence, which means
(contrary to general relativity) they are small perturbations of certain
vacuum solutions. 
These theories are consistent, but phenomenologically uninteresting, as 
in reality one does not observe ten flat dimensions. 
To aleviate this, one can compactify
six of the dimensions, and making the size of the compact factor small
enough, it will be invisible to us yielding a phenomenologically more
acceptable theory. The space time in which the world sheet
resides is then $M^{1,3} \times X$ for some compact manifold $X$. In
order to assure superconformal symmetry here also, it is required that $X$
be Calabi-Yau\footnote{the super- invariance implies K\"ahler, the
  conformal invariance implies Ricci flatness}. 
Such a manifold has two types of moduli: complex structures
and K\"ahler forms. These moduli turn out to be moduli of the physical
theory also, that is, any $X$ (whatever its moduli) gives a consistent
compactification, and in this way, moduli spaces enter also in the
physical theory. These correspond in a sense to ``flat directions of the
potential'', in other words change only the explicit form of the
Lagrangian, not the equations of motion. The ``physical moduli space'' will
in general consist of these geometric moduli in addition to others. 
In order to better understand
these compactifcations, one often compactifies fewer than six dimensions,
an example of which we discuss next. 

\subsection{Toy model: IIA(K3) $\longleftrightarrow$ Het$_{E_8\times
    E_8}$(T$^4$)} 

The notation of the title is meant to indicate that one starts with one of
the five theories above and compactifies on the Calabi-Yau manifold in
parenthesis, in this case of dimension four, leaving a six-dimensional
Minkowski space left as ``space-time''. The arrow indicates the so-called
{\it duality}; this means that the {\it non-perturbative} theory, of which
the supersting theory is an approximation (perhaps it would be better to
say the superstring theory is a certain {\it limit}), is {\it the same} on both
sides. Think of this as meaning there is an underlying theory, with certain
moduli, and at certain moduli points one can approximate the theory by the
two different superstring theories. It is quite involved to list what this
implies physically, but one thing is certain: if this is the case, then
both sides must necessarily have the same moduli space, and this can be
verified and even understood. 

The moduli 
space in question is given by 
\begin{equation}\label{modsp} 
{\cal M} = {\rm SO}(4,20;\fZ) {\large \backslash} {\rm SO}(4,20) {\large /}
{\rm SO(4)}\times {\rm SO(20)}
\end{equation}   
and the heterotic/type II duality reduces in this context to 
giving two different interpretations to this moduli space. The group 
${\rm SO}(4,20;\fZ)$ is the discrete group preserving a particular lattice
in $\fR^{24}$, which is described below. 

The symmetric space on the right-hand side of (\ref{modsp}) 
can be described as the space 
of 4-dimensional subspaces $V$ of the 
space $\fR^{(4,20)}$, by which we mean $\fR^{24}$ 
endowed with a metric of signature (4,20),
on which the metric is positiv definite. 
The space $\fR^{(4,20)}$ contains a (unique) selfdual even integral lattice,  
which we denote by $\gG^{(4,20)} \subset \fR^{(4,20)}$, and which 
has the same signature (in the above notations, $\gG^{(4,20)} \cong
U^{\oplus 4}\oplus E_8^{\oplus 2}$). 
The moduli space ${\cal M}$ then can be described 
as the space of all $V$ up to automorphisms of the lattice 
$\gG^{(4,20)}$. 
 
Now, in the heterotic string interpretation of this space the subspace 
$V$ and its orthogonal complement $V^{\perp}$ are associated to 
gauge fields 
describing the Yang-Mills structure of the string. Generically, the
non-abelian 
gauge group $E_8\times E_8$ of the heterotic string is broken to the
abelian group $U(1)^{16}$ (see for example \cite{Aspin}, section 3.5 and
page 468),
i.e., none of the non-abelian gauge fields persist upon compactification. Of 
special significance are those points in the moduli space 
for which $V^{\perp}$ contains nonzero points of the lattice 
$\gG^{(4,20)}$. If $V^{\perp}$ contains such nonzero points 
for which the normalized length $d_P = -\frac{1}{2}<P,P>$ 
equals unity then physical considerations lead to an enhancement 
of the rank of the Yang-Mills group beyond the rank encountered 
at a generic point of the moduli space. So what is varying here in the
moduli space is the surviving non-abelian gauge group. 

The second, type II string, interpretation of the moduli space 
(\ref{modsp}) is obtained by viewing $\fR^{(4,20)}$ as the 
real cohomology ${\rm H}^*$ of a K3 surface $S$. It is known that type IIA
compactifed on a K3 surface is equivalent to {\it a superconformal
  non-linear sigma models on that K3 surface}. Its moduli space can be
described as follows, see \cite{Aspin}
\S 3 for details on this aspect. From the decomposition
\[ { SO(4,20) \over S(O(4) \times O(20))} \cong {SO(3,19) \over
  S(O(3)\times O(19))} \times \fR^{22} \times \fR_+,\]
this has the following interpretation: the first factor is the space of
{\it Einstein} metrics on a K3 surface, the second factor is the moduli
space of the so-called $B$-field and the final factor is the size. 
It remains to determine the meaning of the classes of norm $-2$ occuring
above. But it is well-known how this occurs. Referring back to section
3.4.2 (and letting our K3 surface be denoted by $X$ as there), we have the
Picard lattice $S_X$ and the transcendental lattice $T_X$. The Picard
lattice is spanned by the K\"ahler class and by rational $-2$ curves, 
which correspond to a generic K3 surface aquiring ordinary double points
and then resolving these to yield a smooth K3 surface. In this respect it
is important to recall that the number of moduli preserving such a Picard
lattice is $20-\rho$, see for example \cite{book}, Proposition 2.3.2. For a
generic {\it algebraic} surface, there are 19 moduli and the Picard number
$\rho$ (the rank of the Picard group) is 1, while for a K3 surface with
$\rho >1$, there are only $20-\rho$ moduli of K3 surfaces which are
deformations of the given one and which have isomorphic Picard lattice. 
The relation in this case to the heterotic theory was observed by Witten in 
\cite{ew95}. Namely, if $S$ aquires a  
singularity of any of the A,D, or E types, a configuration of 
rational curves whose dual graph is the Dynkin diagram of one of the groups
of type A, D or E, collapses 
to a point (these are the vanishing cycles). 
Such rational curves have selfintersection $-2$ and thus 
correspond precisely to the vectors encountered above in the 
heterotic description.

\subsection{String dualities in $D=4$}

Next we consider the compactification of Type IIA on a Calabi-Yau threefold
(down to $D=4$) and the compactification of the heterotic string on a
product $S\times T^2$, where $S$ is a K3 surface. Let us think of this in
the complex category and make the simplifying assumption that $S$ is
elliptically fibered (this is an 18-dimensional family, as compared with the
19-dimensional family of K3 surfaces with some fixed polarization) and that
the Calabi-Yau threefold $X$ is fibered by K3 surfaces. Then, we have two
fibrations:
\[ S \lra \fP^1, \quad \quad X\lra \fP^1,\]
and combining the first with $T^2$, we have two fibrations, one of $X$ and
one of $S\times T^2$, both onto $\fP^1$. It is natural to think of the
above duality for a fixed $t\in \fP^1$, and applying the duality {\it
  fiberwise}. There is a physical argument for this, the so-called {\it
  adiabatic limit}. Note however, that if $S$ is one of our K3 surfaces and
$X$ one of our Calabi-Yau threefolds, for which the modulus of the fiber is
{\it constant}, then there is no argument at all necessary: applying the
duality above fiberwise implies that the gauge group of the heterotic
theory is just the group whose weight lattice is the Picard group of the K3
fiber on $X$. 

For these $D=4$ theories, the heterotic side is quite complicated; there
are vector bundles (background fields) on the K3 surface involved, and a
Higgsing process. Nonetheless, our constructions gives us a good idea of
what the gauge group of a heterotic dual theory would be. In \cite{Alz},
there are chains of such dual theories. We give four examples, three of
which are in \cite{Alz}, the other of which is 
considered in \cite{Alz2}. We start
with the K3 surfaces listed as \# 5,7,8 and 9 in section 3.3.2. We then form
the product with the curve $C$ with weights $(2,1,1)$ and $\ell = 18, 24,
20$ and 42 in the respective cases. The image of the twist map is in these
four cases a weighted hypersurface with weights $(1,1,4,12,18),\
(1,1,6,16,24),\ (1,1,8,10,20),\ (1,1,12,28,42)$ of degrees $36,\ 48,\ 40,\
84$. Looking at the table in section 3.3.2, we see that the singular fibers
I$_0^*$,\ IV$^*$,\ III$^*$ and II$^*$ 
of the elliptic fibrations of the K3 surfaces (the fibers of the
fibration) have dual graphs which are the extended Dynkin diagrams of the
types $D_4,\ E_6,\ E_7$ and $E_8$, respectively, yielding as gauge
groups of the heterotic duals $SO(8),\ E_6,\ E_7$ and $E_8$, respectively. 
For the three cases which occur in \cite{Alz}, these are (up to abelian
factors) indeed the gauge groups of the heterotic dual theories, and for
the remaining case this is verified in \cite{Alz2}, p.\ 133. 

As explained in \cite{Asp2}, the moduli of the two physical theories are
described in more detail as follows. The type IIA string compactified on a
Calabi-Yau threefold $X$ has the following moduli:
\begin{enumerate} \item The dilaton (which governs the string coupling
  constant) and the axion which together form a complex scalar $\Phi_{\rm
  IIA}$;
  \item A metric which is determined by the complex structure of $X$ and
    the cohomology class of a K\"ahler form;
  \item A skew field $B \in {\rm H}^2(X,\fR)/{\rm H}^2(X,\fZ)$;
  \item So-called Ramond-Ramond fields $R\in {\rm H}^{odd}(X,\fR)/{\rm H}^{odd}(X,\fZ)$.
  \end{enumerate}
These moduli split into two types, the so-called {\it vector} moduli and
the {\it hypermultiplet} moduli, as follows:
\begin{itemize}\item[V] The K\"ahler form and the $B$-field together form
  the well-known ``complexified K\"ahler form'' used in mirror symmetry;
  these moduli together form a moduli space $\cM_V$ which is a special
  K\"ahler manifold.  
  \item[H] The complex scalar $\Phi_{\rm IIA}$, the Ramond-Ramond fields
    $R$ and the complex structure of $X$ form the moduli space $\cM_H$ of
    hypermultiplet moduli; this is a quaternionic K\"ahler manifold.
  \end{itemize}
The above descriptions are not valid in the complete quantum theory, but
rather only for certain approximations: $\cM_V$ is valid only in the
``large radius limit'' of $X$, $\cM_H$ is valid only near the
weakly-coupled limit $\Phi_{\rm IIA}\lra -\infty$. 

The heterotic string compactified on a product of a K3 surface $S_H$ and an
elliptic curve $E_H$, $S_H\times E_H$ (which is K\"ahler Ricci-flat) has
the following moduli:
\begin{enumerate} \item The dilaton (which governs the string coupling) and
  the axion, which together form a complex scalar $\Phi_{\rm Het}$.
  \item A Ricci-flat metric on the product $S_H\times E_H$.
  \item A skew-field $B\in {\rm H}^2(S_H\times E_H,\fR)/{\rm H}^2(S_H\times
    E_H,\fZ)$. 
  \item A $G$-bundle on $S_H\times E_H$ with a connection satisfying the
    Yang-Mills equations, where $G$ is the (unbroken) gauge group of the
    heterotic string, either ${\rm Spin}(32)/\fZ_2$ or $E_8\times E_8$.
  \end{enumerate}
Once again, these moduli split into two types, vector and
hypermultiplet. Assume that the $G$ bundle is the product of a
$G_S$-bundle over $S_H$ and a $G_E$-bundle over $E_H$, where $G_S\times
G_E\subset G$ is a subgroup. Then these types can be described as follows:
\begin{itemize}\item[V] The scalar $\Phi_{\rm Het}$, the moduli of the
  $G_E$-bundle over $E_H$, the metric on $E_H$ and the $B$-field on $E_H$
  form the vector muliplet space $\scM_V$. 
  \item[H] The moduli of the $G_S$-bundle on $S_H$, the metric on $S_H$ and
    the $B$-field on $S_H$ form the hypermultiplet moduli space $\scM_H$. 
  \end{itemize}
Again, these descriptions are only valid in certain limits: $\scM_V$ when
$\Phi_{\rm Het} \lra -\infty$ and the area of $E_H$ is large, and $\scM_H$
when the volume of $S_H$ is large. Duality here means essentially matching
these moduli spaces in the two cases. The match of vector moduli is
described above, and the match of hypermultiplet moduli leads to quite
interesting mathematical constructs, for example intermediate Jacobians,
Prym varieties and Deligne cohomology, see \cite{CD}. 

Perhaps the most fascinating aspect of this is determining the K3 surface
explicitly in terms of the Calabi-Yau threefold $X$; for this one considers
a stable degeneration of $X$ into the union of two generalized
Fano-threefolds $X_1 \cup X_2$, and the K3 surface for the heterotic
compactification is the intersection $X_1\cap X_2$.


\subsection{Conifold transitions}
A second problem which is illuminated by our construction is 
the heterotic structure of the so-called conifold transition 
between Calabi-Yau manifolds. Such transitions  are given by the following
construction. Allow a smooth Calabi-Yau threefold $X_t$ depending on a
parameter $t$ to aquire a certain
number of ordinary double points at $t=0$; 
let $X^*(=X_0)$ denote the singular space. Each
of the ordinary double points can be resolved by a small resolution; we
assume that there is at least one of these which is projective (i.e.,
K\"ahler), and let $X^s$ denote such a resolution. 
Schematically we have the following situation:
\[ X_t \llra X^* \llra X^s.\]
The Hodge numbers of $X_t \ (t\neq 0)$
and $X^s$ are related as follows:
\[h^{2,1}(X^s) = h^{2,1}(X_t)-(P-R),\quad \quad h^{1,1}(X^s) = h^{1,1}(X_t)
+R,\]
where $P$ denotes the number of nodes and $R$ denotes the number of
relations between the corresponding vanishing cycles. 
Although this transition passes through a singular space $X^*$, Strominger
has shown \cite{strom} that the physics remains smooth. What happens is
something quite similar to what occurs in the work of Seiberg and Witten: a
massive particle (in this case a black hole) gets massless at the moduli
point of $X^*$ (in the physical theory one always makes a low-energy
approximation, and all massive particles are so heavy that they do not
influence the physics; accordingly they are integrated out of the
Lagrangian), and passing to $X^s$ amounts to the new theory with an
additional massless particle.  
\subsubsection{Splittings}
One way of describing such transitions is by means of {\it splittings},
which are described in terms of complete interesections in products of
weighted projective spaces. As an example of this, consider first a
transversal weighted hypersurface $\fP_{(k_1,k_1,k_2,k_3,k_4)}[d]$, where
$d= 2k_1+k_2+k_3+k_4$, and consider the following threefold in the product 
$\fP_{(1,1)}\times \fP_{(k_1,k_1,k_2,k_3,k_4)}$: 
\[ X_0:= \left\{ \begin{array}{c} p_1(u,y)=u_0Q(y)+u_1R(y)=0 \\ p_2(u,y) =
    u_0S(y)+u_1T(y)=0 
  \end{array} \right\} \inn \fP_{(1,1)}\times
\fP_{(k_1,k_1,k_2,k_3,k_4)}.\]
Schematically this is abbreviated with the following notation: \[X_0 \in
\begin{array}{l} \fP_{(1,1)} \\ \fP_{(k_1,k_1,k_2,k_3,k_4)} 
\end{array} \left[ \begin{array}{cc} 1 & 1 \\ a\cdot k_1 & d-a\cdot k_1
  \end{array} \right],\]
where $a\cdot k_1=deg(Q)=deg(R)$ and $d-a\cdot k_1 = deg(S)=deg(T)$. On the
other hand, consider the determinental variety 
\[ X^* := \{ Q(y)T(y) - R(y)S(y) =0 \} \inn \fP_{(k_1,k_1,k_2,k_3,k_4)}.\]
Clearly $X^*$ is singular for generic choices of $Q,\ R,\ S$ and $T$ when
$Q=R=S=T=0$, which means that $X^*$ generically has isolated singularities,
which one can check are ordinary double points. Furthermore, mapping
$\fP_{(k_1,k_1,k_2,k_3,k_4)}$ into the product $\fP_{(1,1)}\times
\fP_{(k_1,k_1,k_2,k_3,k_4)}$ in the obvious way, it is clear that $X^*$
maps to $X_0$: write the equation defining $X_0$ as \[P(u,y) =
(u_0,u_1)\left( \begin{array}{cc} Q(y) & R(y) \\ S(y) & T(y)
\end{array}\right) = (u_0,u_1)\Pi(y) =0,\]
with a $2\times 2$ matrix $\Pi$. Then $y\in X^* \iff \det(\Pi(y))=0 \iff
\exists_{(u_0,u_1)}$ with $(u_0,u_1)$ is in the kernel of $\Pi(y)$ $\iff$
$P(u,y)=0 \iff (u,y)\in X_0$. The singular $X^*$ can be deformed by adding
some multiple of a transversal polynomial, i.e., by setting 
\[ X_t := \{ t_0( Q(y)T(y) - R(y)S(y) ) + t_1 p_{trans}(y) =0 \} \inn
\fP_{(k_1,k_1,k_2,k_3,k_4)}\quad (t=(t_0,t_1)\in \fP^1) ,\]
and we are in the situation mentioned above: the smooth $X_t$ aquires
singularities of the desired type (ordinary double points) at $t_1=0$, and
this singular $X_0$ can be given a small resolution $X^s$. Once again, we
{\it schematically} describe this process by the symbols
\[ \fP_{(k_1,k_1,k_2,k_3,k_4)}[d] \llra  \begin{array}{l} \fP_{(1,1)} \\
  \fP_{(k_1,k_1,k_2,k_3,k_4)}  
\end{array} \left[ \begin{array}{cc} 1 & 1 \\ a\cdot k_1 & d-a\cdot k_1
  \end{array} \right].\]

\subsubsection{K3 fibrations}
It was shown in \cite{ls95} that such transitions can be  
constructed between K3-fibered Calabi-Yau manifolds, an example 
being provided by the transition
\begin{equation}\label{hetsplit}
\fP_{(1,1,2,4,4)}[12]^{(5,101)} ~\longleftrightarrow ~
\begin{array}{l}\fP_{(1,1)}\\  \fP_{(4,4,1,1,2)}
\end{array}
\left[\begin{array}{cc} 1 & 1 \\ 4 & 8 \\ 
  \end{array} \right]^{(6,70)},
\end{equation}
where the notation on the right denotes a complete intersection manifold 
of codimension two defined by two polynomials of bi-degree 
(1,4) and (1,8) respectively and the superscripts indicate the 
Hodge numbers $(h^{(1,1)},h^{(2,1)})$. 
Note that the smooth hypersurface on the left hand side is the first
example in Table 3, so that the Fermat hypersurface is in fact the image of
an appropriate twist map. In this case we have $\deg(Q)=\deg(R)=4$ and
$\deg(S)= \deg(T)=8$, hence $Q=R=S=T$ has (after rescalings) 32
solutions. Hence the singular $X^*$ above has $P=32$ and $R=1$, which
explains the change in Hodge numbers explicitly from this point of view. 
To see the K3-fibration on the right hand side, one considers the sections
$\gl_0 y_0-\gl_1y_1=0$; an easy calculation shows that these are K3 surfaces
in the family of complete intersections in $\fP_{(1,1)}\times
\fP_{(1,1,2,2)}$ of degrees $(1,2)$ and $(1,4)$.  

\subsubsection{A generalized twist map}
To see that the right hand side above can also be realized as a {\it
  constant modulus} fibration, we can generalize the twist map to this
situation. Define the following rational map:
\begin{eqnarray}\label{gentwist} \Phi: \fP_{(w_0,\ldots,w_n)} \times
  \fP_{(1,1)} \times \fP_{(v_0,\ldots,v_m)} & \lra & \fP_{(1,1)} \times
  \fP_{(w_1v_0,\ldots, w_nv_0,w_0v_1,\ldots, w_0v_m)} \\
  ((x_0,\ldots, x_n),(u_0,u_1),(y_0,\ldots, y_m)) & \mapsto &
  ((u_0,u_1),(y_0^{w_1/w_0}x_1,\ldots,
  y_0^{w_n/w_0}x_n,x_0^{v_1/v_0}y_1,\ldots, x_0^{v_m/v_0}y_m)).\nonumber 
\end{eqnarray}
Let the subvarieties $V_1,\ V_2$ be defined as follows:
  $V_1=\{x_0^{\mu}+p(x_1,\ldots, x_n) =0 \} \inn 
  \fP_{(w_0,\ldots, w_n)};$ 
\begin{equation}\label{vees} \quad V_2 = \left\{ \begin{array}{l}
  p_1({\bf u},{\bf y})= u_0y_1 + u_1\cdot p_{11}(y_1,\ldots, y_m) =0 \\
  p_2({\bf u},{\bf y})= u_0(y_0^{\mu}+p_{20}(y_1,\ldots,y_m))+
  u_1y_1^{\nu-1} =0 
\end{array} \right. \inn \fP_{(1,1)}\times \fP_{(v_0,\ldots, v_m)}.
\end{equation}
This complete intersection has bidegrees $\left[\begin{array}{cc}1 & 1 \\
    v_1 & d-v_1
\end{array} \right]$, 
where $d=\sum_0^m v_i$; $p_{11}$ has degree $v_1$, and
the degree of $p_{20}$ is $\deg(p_{20}) = (\deg( y_1))(\nu-1) = v_1(\nu-1) =
\deg (y_0) \cdot \mu = v_0\cdot \mu = d-v_1$, hence we have the relation
among the various weights:
\[ \mu = {v_1(\nu-1)\over v_0} = {d-v_1 \over v_0}.\] 
Clearly
$V_1$ is invariant under the obvious action of $\fZ_{\mu}$. We claim that
$V_2$ is invariant under the following action of $\fZ_{\nu}$:
\[ ({\bf u,y}) \mapsto ((\eta^{\nu-1}u_0,u_1),(\eta y_0,y_1, \ldots,
y_m))\]
for a generator $\eta$ of $\fZ_{\nu}$. Indeed, the first equation defining
$V_2$ is invariant, while the second gets multiplied by a factor of
$\eta^{\nu-1}$, hence the zero locus is invariant. Set $\ell := {\rm
  gcd}(\mu,\nu)$, then we get an action of $\fZ_{\ell}$ on the product
space as follows:
\[ (\zeta,{\bf x,u,y}) \mapsto ((\zeta^{\nu} x_0,x_1,\ldots,
x_n),(\zeta^{\mu\cdot (\nu-1)} u_0, u_1),(\zeta^{\mu} y_0, y_1, \ldots,
y_m)).\]

Finally, let $X$ be defined in the product \[ \fP_{(1,1)} \times
\fP_{(w_1v_0,\ldots, w_nv_0,w_0v_1,\ldots, w_0v_m)} \hbox{ with coordinates
} ((u_0,u_1),(z_1,\ldots, z_n, t_1,\ldots, t_m))\] by 
\begin{equation}\label{eks} X = \left\{ \begin{array}{l} u_0\cdot t_1
  +u_1p_{11}(t_1,\ldots,t_m) =0 \\
  u_0(p_{20}(t_1,\ldots, t_m)-p(z_1,\ldots, z_n)) ) + u_1\cdot t_1^{\nu-1} =0 
\end{array} \right. . 
\end{equation}
Then a calculation shows that $\Phi(V_1\times V_2) \inn X$,
and that this displays $X$ rationally as a quotient of $V_1\times V_2$ by a
$\boldmu_{\ell}$-operation. Indeed, the first equation defining $X$, in
terms of the coordinates $\bf x,y$, is 
\[ u_0 \cdot x_0^{v_1/v_0}\cdot y_1 + u_1\cdot (x_0^{1/v_0})^{deg(p_{11})}
\cdot p_{11}(y_1,\ldots,y_m) = x_0^{v_1/v_0}(u_0 y_1 + u_1
p_{11}(y_1,\ldots, y_m)),\]
which clearly vanishes for ${\bf (u,y)} \in V_2$. A similar
calculation for the second equation is 
\[ u_0\left( \left(x_0^{1/v_0}\right)^{v_0\cdot \mu} \cdot
  p_{20}(y_1,\ldots, y_m) - \left( y_0^{1/w_0}\right)^{\mu \cdot w_0} \cdot
  p(x_1,\ldots, x_n) \right) + u_1\left( x_0^{v_1/v_0}\right)^{\nu-1} \cdot
  y_1^{\nu-1} ,\]
which again clearly vanishes for ${\bf (x,u,y)} \in V_1\times V_2$.  

\subsubsection{Gauge groups}
Let us apply this to determine the gauge groups on both sides of
(\ref{hetsplit}). For convenience we stick to Fermat polynomials. On the left
hand side we have a constant modulus K3-fibration with fiber the K3 surface
occuring in the first line of the table in section 3.3.2, with six singular
fibers of type IV. Since each of these corresponds to a $A_2$, we get for
the Picard group of the K3 a lattice of type $A_2^6 \oplus H$, where $H$ is
a sum of hyperbolic and abelian factors. As we have already explained,
since the Calabi-Yau threefold is the quotient of the product by a group,
and the action of the group on the K3-surface is of non-Nikulin type, it
preserves the Picard lattice, hence this is, up to abelian factors, the
gauge group of the theory compactified on that Calabi-Yau.
To do the same for the right hand
side, we must first study the K3 fiber, which is the complete interesection
of type $(1,2)$, $(1,4)$ in $\fP_{(1,1)}\times \fP_{(1,1,2,2)}$. Recall
that we can determine the Picard group from the elliptic fibration; in this
case we first describe this K3 as the image of a twist map. This is given
by the above constuction, setting $(v_0,v_1,v_2)=(1,1,1),\ (w_0,w_1,w_2)
=(2,1,1)$, $p(x_1,x_2)=x_1^4+x_2^4$, $p_{11}(y_1,y_2)=y_2$ and
$p_{20}(y_1,y_2) = y_1^2$. In other words, let $V_1 = \{
x_0^2+x_1^4+x_2^4=0\}\inn \fP_{(2,1,1)}$ (so $\mu =2$),
and let $V_2$ be the elliptic curve defined as follows:
\[ V_2 = \left\{ \begin{array}{l} u_0\cdot y_1+u_1\cdot y_2 =0 \\ 
  u_0(y_0^2+y_1^2) + u_1\cdot y_1^2 =0 
  \end{array} \right. \inn \fP_{(1,1)}\times \fP_{(1,1,1)}.\]
The image of $V_1\times V_2$ under the twist map is the complete
intersection $W$ defined as follows
\[ W = \left\{ \begin{array}{l} u_0t_1+u_1t_2 =0\\
    u_0(t_1^2-(z_1^4+z_2^4))+u_1\cdot t_1^2 =0
  \end{array} \right. \inn \fP_{(1,1)}\times \fP_{(1,1,2,2)}.\]
In this case, $\nu=3$ and $\ell =6$, so we have a group of order six acting
on the product, with a subgroup of order 3 acting on the elliptic
curve. The singular fibers of the fibration are hence of types {\bf IV} or
{\bf IV$^*$}. It is easy to see that at the four zeros of $p(x_1,x_2)$, we
have singular fibers of type ${\bf IV}$. The reader may check that there
is, in addition, a singular fiber of type ${\bf IV^*}$. Hence, in this case
the Picard group is $A_2^4\oplus E_6$. 

As we have described above, this is the same as a determinental
hypersurface in $\fP_{(1,1,2,2)}$, which aquires $\deg(Q)\cdot \deg(R)\cdot
\deg(S)\cdot \deg(T) = 4$ ordinary double points. However, these four
ordinary double points all lie on certain rational curves, and after
resolution of singularities, they are components of the $\bf IV^*$ fiber. 
Since we
have a constant modulus fibration, fiberwise duality implies that
$A_2^4\oplus E_6$ is
indeed the gauge group of the heterotic string {\it after} the conifold
transition. This kind of information is new and exciting from the point of
view of physics.

\subsection{Elliptic fibrations} 
More recently not only K3-fibrations have become of importance 
in string theory, but also elliptic fibrations in various dimensions.
In the framework fo F-theory \cite{MV} a more general type of 
compactification of the ten-dimensional string (of type IIB) is 
considered in which the dilaton field of the string is not 
assumed to be constant, as in conventional compactifications, 
 but instead can vary. 
In this solution of the string equations this dilaton is assumed 
to take values in an elliptic curve, thus leading to description 
of these vacua as 12-dimensional elliptic fibrations. 
Interesting theories then can be described by elliptically fibered 
Calabi-Yau manifolds of complex dimensions two, three and four. 

These new vacua of the type IIB string are conjectured to be dual 
to the most general compactification possible for the heterotic 
string, based on stable vector bundles with vanishing first 
Chern class. 
Of particular interest in this context 
 are Calabi-Yau fourfolds which not only 
are elliptic but also K3-fibered. Such spaces lead to a four-dimensional 
compactification of F-theories which are expected to be dual 
to stable vector bundles over Calabi-Yau threefolds which are elliptic. 
 Clearly our twist construction provides a systematic method for 
building and analyzing such varieties. 

In this general framework the question of phase transitions rises 
again in the context of the possible connectedness of these 
fourdimensional string theory ground states. This problem has 
been addressed in \cite{bls97} in the framework of 4D F-theory. 
It was shown there that there indeed exists a generalization 
of the conifold transition for Calabi-Yau fourfolds. In the 
context of these higher dimensional varieties however the 
singularities are no longer nodes but the 
manifolds degenerate at curves of in general high genus. 

A simple example of such a transition between CY$_3$ fibered 
fourfolds is provided by \cite{bls97} 
\begin{equation}
\fP_{(8,8,4,2,1,1)}[24] \llra  
\matrix{\fP_{(1,1)}\hfill \cr \fP_{(8,8,4,2,1,1)}\cr }
\left[\matrix{1&1\cr 8&16\cr}\right].  
\end{equation} 
 In this transition the singular locus is given by the smooth 
curve $\gS =\fP_{(4,2,1,1)}[16~~16]$ of genus $g=385$. 
At the singular configuration of this space the CY$_3$ fiber 
degenerates into a conifold configuration with 32 nodes. 
We can obtain this singular fiber by either deforming the 
generic hypersurface fiber $\fP_{(4,4,2,1,1)}[12]$ which is 
obtained from the twist map 
\begin{equation}
\fP_{(2,1,1)}[12]\times \fP_{(2,2,1,1)}[6] ~\lra ~ 
\fP_{(4,4,2,1,1)}[12]
\end{equation}
i.e. by collapsing 32 three-cycles, or by collapsing 32 two-cycles 
of the generic quasismooth codimension two complete intersection 
CY which leads to the codim$_{\fC}=2$ fourfold via the 
twist map 
\begin{equation}
\fP_{(2,1,1)}[8] \times
  \matrix{\fP_{(1,1)}\hfill \cr \fP_{(4,4,2,1,1)}\cr}
  \left[\matrix{1&1\cr 4&8\cr}\right]
~\lra ~
 \matrix{\fP_{(1,1)}\hfill \cr \fP_{(8,8,4,2,1,1)}\cr}
       \left[\matrix{1&1\cr 8&16\cr}\right]
     \end{equation}
     
From the analysis above we expect that via our twist map 
many of the aspects of the conifold transitions have 
important implications for the transitions among Calabi-Yau 
fourfolds and, via the conjectured duality between F-theory 
and the heterotic string compactified on stable vector bundles 
$V~\rightarrow ~$CY$_3$ also for transitions between stable 
vector bundles.

\smallskip

\noindent Authors' addresses:

\noindent
{\it Bruce Hunt \\ Max-Planck-Institut f\"ur Mathematik in den
   Naturwissenschaften \\
   Inselstr. 22-26, D-04103 Leipzig \\ \\ Rolf Schimmrigk \\ Physics
   Institute, Bonn University \\  Nu\ss allee 12, D-53115 Bonn \\
   and \\ Dept.\ of Physics and Astronomy, Indiana University South Bend
   \\ 1700 Mishawaka Av., South Bend, IN 46634, USA }

\end{document}